\documentclass[reqno,11pt]{amsart}
\usepackage{amsmath, amssymb, amsthm,amsfonts} 
\usepackage[english]{babel}
\usepackage{bbm}
\usepackage{graphicx}
\usepackage{url}
\usepackage{epstopdf}
\usepackage[a4paper,bindingoffset=0.5cm,left=2cm,right=2cm,top=2.5cm,bottom=2cm,footskip=.8cm]{geometry}

\usepackage{hyperref}
\hypersetup{
    colorlinks=true,
    linkcolor=blue,
    filecolor=blue,      
    urlcolor=blue,
}

\usepackage{tikz}
\usetikzlibrary{arrows, automata,positioning,calc,shapes,decorations.pathreplacing,decorations.markings,shapes.misc,petri,topaths}
  \usepackage{pgfplots}
  \pgfplotsset{compat=newest}
  \usetikzlibrary{plotmarks}
  \usepackage{grffile}
\newlength\figureheight
  \newlength\figurewidth
\pgfplotsset{%
    tick label style={font=\scriptsize},
    label style={font=\footnotesize},
    legend style={font=\footnotesize},
         every axis plot/.append style={very thick}
}
\usetikzlibrary{math,shapes.geometric}

\usepackage{rotating}
\usepackage{amsbsy,enumerate}
\usepackage{graphicx}
\usepackage{ccaption}
\usepackage{comment}
\usepackage{mathrsfs} 
\usepackage{diagbox}
\usepackage{floatrow}

\newcommand{\bs}{\boldsymbol}

\newcommand{\vb}{\vspace{3.2mm}}

\newcommand{\tred}[1]{{\textcolor{black}{#1}}}
\newcommand{\tredd}[1]{{\textcolor{black}{#1}}}

\allowdisplaybreaks

\newtheorem{lemma}{Lemma}

\newtheorem{theorem}{Theorem}
\newtheorem{remark}{Remark}
\newtheorem{example}{Example}

\newtheorem{proposition}{Proposition}

\makeatletter
\renewcommand{\fnum@figure}[1]{\textbf{\figurename~\thefigure}. }
\renewcommand{\fnum@table}[1]{\textbf{\tablename~\thetable}. }
\makeatother

\begin{document}

\title[Adaptive Scheduling in Service Systems: A Dynamic Programming Approach]{Adaptive Scheduling in Service Systems:\\ A Dynamic Programming Approach}

\author{Roshan Mahes, Michel Mandjes, Marko Boon, and Peter Taylor}

\begin{abstract} This paper considers appointment scheduling in a setting in which at every client arrival the schedule of all future clients can be adapted. Starting our analysis with an explicit treatment of the case of exponentially distributed service times, we then develop a phase-type-based approach to also cover cases in which the service times' squared coefficient of variation differs from 1. The approach relies on dynamic programming, with the state information being the number of clients waiting, the elapsed service time of the client in service, and the number of clients still to be scheduled. The use of dynamic schedules is illustrated through a set of numerical experiments, showing (i)~the effect of wrongly assuming exponentially distributed service times, and (ii)~the gains (over static schedules, that is) achieved by rescheduling. 

\vb


\vb

\noindent
{\sc Affiliations.} 
\noindent
{\it Roshan Mahes} and {\it Michel Mandjes} are with the Korteweg-de Vries Institute for Mathematics, University of Amsterdam, Science Park 904, 1098 XH Amsterdam, The Netherlands, and Amsterdam Business School, Faculty of Economics and Business, University of Amsterdam, Amsterdam, the Netherlands. MM is also with E{\sc urandom}, Eindhoven University of Technology, Eindhoven, The Netherlands; his research is partly funded by the NWO Gravitation project N{\sc etworks}, grant number 024.002.003. 

\noindent
{\it Marko Boon} is with the Department of Mathematics and Computer Science, Eindhoven University of Technology, P.O. Box 513, 5600 MB Eindhoven, The Netherlands.

\noindent
{\it Peter Taylor} is with the School of Mathematics and Statistics, University of Melbourne, Parkville, Victoria 3010,
Australia. PT's research is partly funded by the ARC through Laureate Fellowship FL130100039 and the ARC Centre of Excellence for the Mathematical and Statistical Frontiers (ACEMS).

\vb

\noindent
{\sc Acknowledgments.} The authors would like to thank John Gilbertson (for discussions on the homogeneous exponential case and making \cite{GIL} available), Alex Kuiper (for providing Figure \ref{fig:alex}) and Ruben Brokkelkamp (for numerical support). \tred{In addition we thank the anonymous reviewers and associate editor for helpful and constructive comments.}

\end{abstract}

\maketitle

\newpage

\section{Introduction}
In a broad range of service systems, the demand is regulated by working with scheduled appointments. 
When setting up such a schedule, the main objective is to properly weigh the interests of the service provider and its clients. One wishes to efficiently run the system, but at the same time,  the clients should be provided with a sufficiently high level of service. This service level is usually phrased  in terms of the individual clients' waiting times, whereas the system efficiency is quantified by the service provider's idle time. In order to generate an optimal schedule, one is faced with the task of identifying the clients' arrival times that minimize a cost function that encompasses the expected idle time as well as the expected waiting times. 

The above framework has been widely adopted in  healthcare \cite{CAY, GUP}, but there is ample scope for applications in other domains as well. When thinking for instance of a delivery service, clients are given times (or time windows) at which the deliverer will hand over the commodity (parcel, meal, etc.). In this case, idle times manifest themselves if the deliverer happens to arrive earlier than scheduled, whereas clients will experience waiting times when the deliverer arrives later than scheduled. Other examples in which the above scheduling framework can be used can be found in applications in logistics, such as the scheduling of ship docking times \cite{W2}.

In the literature, the dominant approach considers {\it static} schedules: the clients' arrival times (chosen to minimize the given cost function) are determined a priori, and are not updated on the fly; examples of such studies are, besides the seminal paper \cite{BAI}, \cite{HM, KUIP, PEGD, W2}. One could imagine, however, that substantial gains can be achieved if one has the flexibility to now and then adapt the schedule. For instance, if one is ahead of schedule or lagging behind, one may want to notify the clients to provide them with a new appointment time, so as to better control the cost function. In the context of a delivery service, the client will appreciate being given updates on the delivery time, so as to better plan her other activities of the day. \tredd{Besides,} the option of updating the schedule enables the deliverer to mitigate the risk of losing time because of arriving at the client's address before the scheduled time. \tredd{Throughout this paper we will phrase the problem in terms of arriving clients, their appointment epochs, and their service times.} 

In the present work, we consider the setting in which the schedule is adapted at any client arrival time. \tredd{In this respect, it
should be noted that, evidently, for any specific application that one would like to analyze, domain-specific
details need to be added.} \tred{The approach that we develop 
is particularly useful in situations in which clients are essentially already on standby, but benefit from being given a precise estimate of their appointment time. A few examples in which our setup could be used are:}
\begin{itemize}
    \item[$\circ$] \tred{Consider a hospital in which a certain type of surgery is performed. The patients are already on standby (i.e., physically present at the hospital, say, residing \tredd{in a day surgery waiting room.} Then with our algorithm, whenever a \tredd{patient's surgery starts}, the start time of the next patient's surgery can be reliably scheduled. The clear benefit is that it is preferred that the patient spends her waiting time in the \tredd{waiting room} rather than in the immediate vicinity of the operating theater. Furthermore, it is often the case that various preparatory tasks need to be performed, such as commencing an anesthetic regime, that can be better timed when there is a more precise estimate of the operating time.}

    \item[$\circ$] \tredd{Consider a healthcare facility, say an operating theater, that is used by multiple doctors, say on-call surgeons. In order to maximally efficiently use the operating theater, the surgeons can be given `appointment times', i.e., times at which they have to be at the facility to potentially start their surgery. In the language of our scheduling problem, the clients are now the surgeons, each service time corresponds to the total time that it takes a specific surgeon to treat a consecutive stream of patients that were assigned to her, waiting time is the time the surgeon has to wait because of a preceding doctor still being busy, and idle time corresponds to time intervals in which the healthcare facility is not used. In the most basic variant, the patients could be assumed to be on standby, \textcolor{black}{while in more advanced variants two stages of waiting are incorporated in the cost function: (1) time spent in the ‘standby mode’ in which the customer has to be at the service location
or close to their home, and (2) the `true' waiting time, in which service is imminent}.
}
    
    \item[$\circ$] \tred{Consider the situation of a parcel delivery service, \tredd{in which couriers drop parcels according to some predetermined optimized route}. Typically, clients are provided a relatively wide window, devised as a single appointment time with a confidence interval, in which they are supposed to be \tredd{close to their home}. With our approach, they can be given a refinement \tredd{at which time they actually need to be at home. For example, this can be done at the moment the delivery person enters the previous zone (e.g., a street, a zip code, or a small neighborhood). The (random) time between the first parcel delivery of two subsequent zones can be interpreted as a `service time' in our model. This comprises both the handling time and travel time of the driver. Sending this real-time update on the delivery enables the client to be aware of when to step outside, thus reducing the handling time of the delivery driver.}}

    \item[$\circ$] \tred{Various other applications can be thought of. Also in the setting of, say, a repairperson or a plumber, clients are typically given wide intervals in which they have to be present, but they would be pleased to be sent (reliable) adjustments on the fly.} \tredd{In this context, a client's `service’ should be thought of as capturing the time until the repairperson is available to start working on the next service (that is, the travel time is included).}
\end{itemize}
\tred{Besides being applicable in the settings mentioned above, we would like to stress that our approach provides a benchmark for what one could win by adaptively updating the appointment schedule. \tredd{Concretely, with adaptive schedules becoming increasingly important, our numerical output can be used to quantify the potential gain. Our machinery also facilitates the comparison of (conceptually or computationally) simpler adaptive scheduling approaches with the `ideal' DP-based one; as such it is useful when exploring the `complexity/efficiency tradeoff'.}}

The approach followed relies on dynamic programming \cite{BELL}: the new schedule is evaluated by recursively reducing it to simpler sub-problems of a lower dimension. Importantly, the schedule update is performed by making use of all information available at the client arrival time under consideration, namely 
(a)~the number of clients who are already waiting at that time, (b)~the elapsed service time of the client in service (if any), and (c)~the number of clients still to be scheduled. 
In many appointment scheduling studies, such as \cite{HM, PEGD}, service times are assumed to be exponentially distributed, having the advantage of facilitating more explicit analysis. \tredd{In particular, due to the memoryless property, it is not necessary to track information about the service time as in (b) above.}

Empirical studies, such as~\cite{CV} for the healthcare setting, however, show that in various practical contexts, this exponentiality assumption is not justified.
In \cite{VUY}, this is resolved by choosing a discrete-time setting that obtains accurate predictions for any known service-time distribution at low computational cost. 
As the exact service-time distribution is generally unknown, this led us to develop a method that can in principle deal with any type of service-time distribution. 
\tredd{
In our approach, the service-time distributions are parameterized by their mean and variance (or, more conveniently, by their mean and squared coefficient of variation (SCV)). 
Then the idea is to generate a phase-type distribution with this mean and 
variance. Informally, phase-type distributions are effectively generalizations of exponential distributions, thus offering considerably more flexibility but at the same time still allowing 
explicit computations. 
}

\tredd{
The main complication of working with non-exponential distributions lies in the fact that the state of the system has to contain some information about the service time of the client in service. 
It could, for example, be the {\it elapsed} service time or the {\it residual} service time. In the case when the non-exponential distribution is of phase type, it suffices to track the current phase of the service time.}
This complication leads to serious computational challenges: due to the large state space comprising (a), (b), and (c), the numerical evaluation of the schedule is far from trivial. 

We proceed by detailing the paper's main contributions. 
\begin{itemize}
    \item[$\circ$]
    In the first place, we set up a framework for dynamic appointment scheduling.  As first steps, we deal with systems with  exponentially distributed service times, for which still relatively explicit computations can be performed. 
    Importantly, however, our approach is then extended to service-time distributions with an arbitrary SCV. It requires various explicit calculations for queues with deterministic interarrival times and service times that correspond to either a mixture of Erlang distributions or a hyperexponential distribution. Our framework revolves around a dynamic-programming recursion.
    The focus is primarily on the homogeneous case, i.e., the situation that the clients' service times are equally distributed; in the exponential case, however, we also deal with the heterogeneous case, as some of the underlying computations help in the analysis of the case of a general SCV.
    \item[$\circ$]
The second main contribution concerns our computational approach. It relies on a delicate combination of  various numerical procedures, some of which are applied in a rather pragmatic manner. These procedures include efficient search algorithms that are required to determine the optimal appointment times for a given instance, storing a sizeable set of dynamic schedules that correspond to specific parameter settings, and various interpolation and extrapolation steps. All procedures that we applied are backed by extensive testing. \item[$\circ$] Our approach is illustrated by a set of numerical experiments. In the first place, we assess to what extent assuming exponentiality leads to performance degradation. In the second place, we quantify the gain of dynamic schedules over  static, precalculated schedules. This is done as a function of the SCV, for various values of the parameter that weighs the cost of waiting relative to the cost of server idleness.   
\end{itemize}
\tredd{
In the setup presented in this paper, an initial schedule is computed before the start of the first appointment. 
This schedule is dynamically updated whenever a new client arrives. It is up to the scheduler to determine how frequently schedule updates should be communicated, as this may depend on the application area, the lengths of the service times and client travel times (if applicable).}
In the healthcare context, for instance, our concept of dynamic schedules is primarily useful when dealing with treatments that take relatively long, as the patient should have enough time to appear at the healthcare facility. For similar reasons, in the parcel delivery context, our setup is well-suited in situations in which the time between subsequent deliveries is sufficiently long. In addition, applications are envisaged in the context of e.g.\ electricians, plumbers, heating engineers, etc., in which a typical service takes a relatively long time. In practical situations, one may schedule the next client in the way described in our paper, but still provide all subsequent clients with an estimate of the time at which their service is expected to start (which is then later adjusted).

Dynamic schedules have been studied extensively, but with a focus on settings that crucially differ from the one described above. Some papers study what could be called `advance admission scheduling', in which a planner is to dynamically assign appointments for future days, taking into account the clients' preferences and future demand; see for example 
\cite{HG, WF, ZAN}. Another related domain concerns the sequencing of the clients: whereas in our framework the order of the clients is fixed, one could also study the question of identifying the best sequence \cite{BER, MAK}. In this respect, there is empirical evidence and formal backing to opt for the rule that sequences clients in order of increasing variance of their service durations \cite{KEMP, KONG}. \tred{A relevant related research line concerns settings in which there are routine \tredd{clients}, who can be assigned an appointment time in advance, in combination with last-minute \tredd{clients} \cite{CR, ED}; this problem has been formulated in terms of a (multistage) stochastic linear program. 
\tredd{Our type of modeling has only recently been introduced in the delivery literature; in the broader context of home delivery services, e.g.\ \cite{TS, ZWW} work with a setup that is very similar to ours, viz.\ an objective function comprising a weighted sum of mean busy and idle times (where it is noted that \cite{ZWW} has a strong focus on the routing problem, i.e., finding the order that minimizes the objective function, while in our study it is the objective to exploit real-time information to reduce costs by updating the schedule).}
Finally, we would like to mention the significant recent contributions \cite{WLW}, focusing on managing walk-in clients, and \cite{ZY}, exploiting the concept of multimodularity to set up a highly general framework e.g.\ no-shows, non-punctuality, and walk-ins.}

This paper is organized as follows. Section \ref{MF} introduces the modeling framework, defining the dynamics of the underlying system, and formally defining the cost function. We then focus on schedules with exponentially distributed service times: the homogeneous case is dealt with in Section \ref{homexp}, whereas the heterogeneous case can be found in Section \ref{hetexp}. In Section \ref{phase} we develop our phase-type-based framework, setting up the dynamic programming routine that deals with distributions with a general SCV. Numerical experiments are presented in Section \ref{num}.  \tred{In Section~\ref{sec:discu} we discuss robustness properties, as well as a computationally attractive variant in which the so-called value-to-go is ignored.}
We conclude the paper with some final remarks.

\section{Modeling framework}
\label{MF} 
In this section, we describe the modeling framework that will be used in this paper. It is a widely accepted setup that is intensively used in the appointment scheduling literature; cf.\, for instance, \cite{HM, KK, KUIP, PEGD, W1, W2}.

We consider a sequence of $n\in{\mathbb N}$ clients with service times that are represented by the independent, non-negative random variables $B_1,\ldots, B_n$. The idea is that we schedule the jobs one by one, in the sense that at the moment client $i$ enters the system, the arrival epoch of client $i+1$ is scheduled. This we do relying on dynamic programming, with a pivotal role being played by the Bellman equation. At (immediately after, that is) each of the arrival epochs the {\it state} of the system is the number of clients waiting and the elapsed service time of the client in service, where one in addition knows the characteristics of the clients that are still to be served. 
Clients are assumed to arrive punctually. As soon as a client's service has been completed, the server proceeds by serving the next client (if present). 

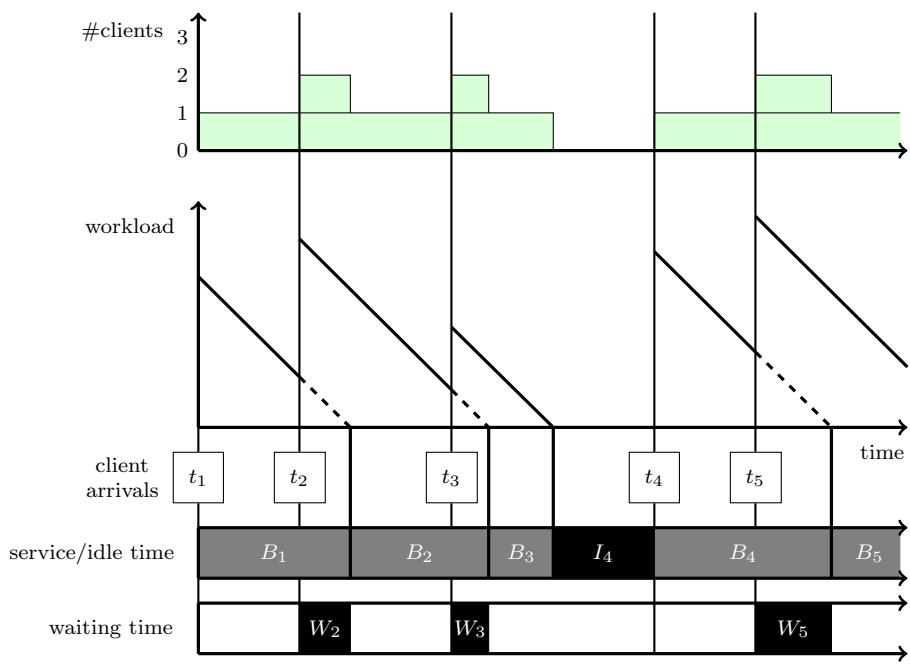
\begin{figure}[h!]
\centering
\begin{tikzpicture}

\node at (-1, 8.25) {\scriptsize \#clients};

\draw[fill, green!15] (0, 6.67) rectangle (4.67, 7.17);
\draw[fill, green!15] (1.33, 7.67) rectangle (2, 7.17);
\draw[fill, green!15] (3.33, 7.67) rectangle (3.82, 7.17);
\draw[fill, green!15] (6, 6.67) rectangle (9.23, 7.17);
\draw[fill, green!15] (7.33, 7.67) rectangle (8.33, 7.17);

\draw (0, 7.17) -- (4.67, 7.17) -- (4.67, 6.67);
\draw (1.33, 7.67) -- (2, 7.67) -- (2, 7.17);
\draw (3.33, 7.67) -- (3.82, 7.67) -- (3.82, 7.17);
\draw (6, 7.17) -- (9.23, 7.17);
\draw (7.33, 7.67) -- (8.33, 7.67) -- (8.33, 7.17);

\draw[thick] (0, 0) -- (0, 5);
\draw[thick] (1.33, 0) -- (1.33, 8.5);
\draw[thick] (3.33, 0) -- (3.33, 8.5);
\draw[thick] (6, 0) -- (6, 8.5);
\draw[thick] (7.33, 0) -- (7.33, 8.5);

\draw[very thick, ->] (0, 6.67) -- (0, 8.5);
\draw [very thick, ->] (0, 6.67) -- (9.33, 6.67);

\node [left] at (0, 6.67) {\scriptsize $0$};
\node [left] at (0, 7.17) {\scriptsize $1$};
\node [left] at (0, 7.67) {\scriptsize $2$};
\node [left] at (0, 8.17) {\scriptsize $3$};


\node at (-0.9, 5.67) {\scriptsize workload};
\node at (9, 2.67) {\scriptsize time};

\draw[very thick, ->] (0, 3) -- (9.33, 3);
\draw[very thick, ->] (0, 3) -- (0, 6);

\draw[very thick] (0, 5) -- (1.33, 3.67); 
\draw[dashed, very thick] (1.33, 3.67) -- (2, 3);
\draw[very thick] (1.33, 5.5) -- (3.33, 3.5); 
\draw[dashed, very thick] (3.33, 3.5) -- (3.82, 3);
\draw[very thick] (3.33, 4.33) -- (4.67, 3); 
\draw[dashed, very thick] (7.33, 4) -- (8.33, 3); 
\draw[very thick] (6, 5.33) -- (7.33, 4);
\draw[very thick] (7.33, 5.8) -- (9.33, 3.8); 


\node at (-1, 2.5) {\scriptsize client};
\node at (-1, 2.17) {\scriptsize arrivals};

\draw [fill, white] (-0.33, 2) rectangle (0.33, 2.67);
\draw (-0.33, 2) rectangle (0.33, 2.67);
\node at (0, 2.33) {\scriptsize$t_1$};
\draw [fill, white] (1, 2) rectangle (1.67, 2.67);
\draw (1, 2) rectangle (1.67, 2.67);
\node at (1.33, 2.33) {\scriptsize$t_2$};
\draw[fill, white] (3, 2) rectangle (3.67, 2.67);
\draw (3, 2) rectangle (3.67, 2.67);
\node at (3.33, 2.33) {\scriptsize$t_3$};
\draw[fill, white] (5.67, 2) rectangle (6.33, 2.67);
\draw (5.67, 2) rectangle (6.33, 2.67);
\node at (6, 2.33) {\scriptsize$t_4$};
\draw[fill, white] (7, 2) rectangle (7.67, 2.67);
\draw (7, 2) rectangle (7.67, 2.67);
\node at (7.33, 2.33) {\scriptsize$t_5$};


\node at (-1.4, 1.33) {\scriptsize service/idle time};

\draw[fill, gray] (0, 1) rectangle (9.23, 1.67);
\draw[fill, black] (4.67, 1) rectangle (6, 1.67);

\node[white] at (1, 1.33) {\scriptsize$B_1$};
\node[white] at (2.91, 1.33) {\scriptsize$B_2$};
\node[white] at (4.25, 1.33) {\scriptsize$B_3$};
\node[white] at (5.33, 1.33) {\scriptsize$I_4$};
\node[white] at (7.17, 1.33) {\scriptsize$B_4$};
\node[white] at (8.83, 1.33) {\scriptsize$B_5$};

\draw[very thick] (0, 1) -- (0, 1.67);
\draw[very thick, ->] (0, 1) -- (9.33, 1);
\draw[very thick, ->] (0, 1.67) -- (9.33, 1.67);

\draw[very thick] (2, 1) -- (2, 3);
\draw[very thick] (3.82, 1) -- (3.82, 3);
\draw[very thick] (4.67, 1) -- (4.67, 3);
\draw[very thick] (8.33, 1) -- (8.33, 3);


\node at (-1.15, 0.33) {\scriptsize waiting time};

\draw [very thick] (0, 0) -- (0, 0.67);
\draw [very thick, ->] (0, 0) -- (9.33, 0);
\draw [very thick, ->] (0, 0.67) -- (9.33, 0.67);

\draw[fill, black] (1.33, 0) rectangle (2, 0.67);
\draw[fill, black] (3.33, 0) rectangle (3.82, 0.67);
\draw[fill, black] (7.33, 0) rectangle (8.33, 0.67);

\node[white] at (1.67, 0.33) {\scriptsize$W_2$};
\node[white] at (3.58, 0.33) {\scriptsize$W_3$};
\node[white] at (7.83, 0.33) {\scriptsize$W_5$};
\end{tikzpicture}
\caption{\label{fig:alex}Illustration of  key quantities in appointment scheduling. Figure courtesy of A.~Kuiper \cite{AThesis}.}

\end{figure}

Our objective function, or cost function, reflects the interests of both the clients and the service provider, in that it captures the `disutilities' experienced by both parties involved. Both components are weighed to reflect their relative importance. The clients' disutility is measured through their waiting times, and the service provider's disutility through the queue's idle times. Concretely, the cost function that we will be working with is
\begin{equation}\label{OBJ}\omega   \sum_{i=1}^n {\mathbb E} \, I_i + (1-\omega)  \sum_{i=1}^n {\mathbb E} \, W_i,\end{equation}
where $I_i$ is the idle time and $W_i$ the waiting time associated with client $i$. The arrival epoch of client $i$ will be denoted by $t_i$, where evidently $t_1$ can be set equal to $0$ (and obviously $I_1=0$); the corresponding interarrival times $t_{i}-t_{i-1}$ are denoted by $x_i$ (where we set $x_1 = 0$). For a pictorial illustration, we refer to Figure \ref{fig:alex}. 

Note that we can write the schedule's {\em makespan} (i.e., the time until the last client has left the system) in multiple ways:
\[\sum_{i=1}^n I_i + \sum_{i=1}^n B_i = t_n+ W_n + B_n=\sum_{i=1}^{n} x_i +W_n+B_n;\]
the left-hand side follows from the observation that every point in time during the makespan belongs to either an idle time or a service time, whereas the other two expressions follow from the fact that the makespan equals the arrival epoch of the last client increased by her waiting time (if any) and her service time.  Realizing that ${\mathbb E}\, B_i$ are given numbers, the consequence of the above identity is that we can equivalently consider the cost function 
\[\omega   \left(\sum_{i=1}^{n}x_i +{\mathbb E}\,W_n\right) + (1-\omega)  \sum_{i=1}^n {\mathbb E} \,W_i.\]

Let $C_i(k,u)$ denote the minimal cost incurred from the arrival of the $i$-th client (until the end of the schedule), given there are $k$ clients in the system immediately after the arrival of the $i$-th client, and the client in service has an elapsed service time $u\geqslant 0$. 
The objective is to minimize \eqref{OBJ}, which is equivalent to evaluating $C_1(1,0)$. In the next section, we start our analysis with the easiest case: the homogeneous exponential case, i.e., each $B_i$ is exponentially distributed with mean $\mu^{-1}$, to gradually proceed to the setup with a general SCV in Section \ref{phase}.

\section{Homogeneous exponential case}\label{homexp}
When the service times are exponentially distributed, the state of the system is just the number of clients waiting. Observe that the elapsed service time of the client in service is irrelevant, by virtue of the memoryless property of the exponential distribution. In this section, we set up an efficient technique to determine the optimal arrival time of the next client in case the $B_i$ stem from the same exponential distribution, for the dynamic scheduling mechanism we defined in Section \ref{MF}. 
\tredd{
The main result of this section is Theorem \ref{t1}, where we present an efficient, dynamic-programming \cite{BELL} based algorithm to dynamically find the optimal arrival times. 
The algorithm dynamically evaluates the cost function  by conditioning on the state of the system (i.e., the number of clients in the system) at  arrival epochs. To evaluate this cost function, we first determine the contribution 
due to the idle and waiting times 
(computed in Proposition~\ref{P1}), and then determine the transition probabilities (in Lemma~\ref{L1}).
}

Suppose we wish to evaluate the cost between the arrival of the $i$-th and $(i+1)$-st client. For ease, we (locally) identify time $0$ with the arrival of this $i$-th client, and we assume that, immediately after this arrival, there are $k$ clients in the system (where $k=1,\ldots,i$). We let $t$ be the time at which the $(i+1)$-st client is scheduled to arrive. Let $N_s$ denote the number of clients in the system at time $s\in[0,t]$, including the client in service.

The first observation is that there are essentially two \tredd{(mutually exclusive) scenarios upon the arrival of client $i+1$: (i)~the newly arrived client immediately gets treatment, i.e., the system idles before $t$, and (ii)~the newly arrived client has to wait in the queue, i.e., there is still a previously arrived client in service at time $t$.} A second observation is, considering the objective function that we defined in \eqref{OBJ}, in the former scenario there is a contribution to the mean idle time as well as the mean waiting time, whereas, in the latter scenario, only the mean waiting time increases. 

\begin{itemize}
\item[$\circ$]
The contribution of the idle time (in scenario (i)) to the cost function, due to the interval $[0,t]$, is $\omega f_k(t)$, with
\begin{align}f_k(t)&:= \int_0^t {\mathbb E}\big(
\mathbbm{1}_{\{N_s=0\}}
\,|\, N_{0+}=k\big) {\rm d}s= \int_0^t {\mathbb P}\big(N_s=0\,|\, N_{0+}=k\big) {\rm d}s.\label{FF1}\end{align}
\item[$\circ$] The contribution of the waiting time (in both scenarios) to the cost function, due to the interval $[0,t]$, is $(1-\omega)g_k(t)$, with
\begin{align}\nonumber g_k(t)&:= \int_0^t\sum_{\ell=0}^{k-1} (k-\ell-1) {\mathbb E}\big(
\mathbbm{1}_{\{N_s=k-\ell\}}
\,|\, N_{0+}=k\big) {\rm d}s\\
&=\int_0^t\sum_{\ell=0}^{k-1} (k-\ell-1) {\mathbb P}\big(N_s=k-\ell\,|\, N_{0+}=k\big) {\rm d}s.\label{GG1}
\end{align}
Here we used that if there are $\ell$ clients in the system, then $\ell-1$ of them are waiting.
\end{itemize}

The quantities $f_k(t)$ and $g_k(t)$ can be evaluated by means of a direct calculation of the integrand appearing in the right-hand sides of \eqref{FF1} and \eqref{GG1}. The main idea is to condition on the number of service completions up to time $s$. To this end, observe that the event $\{N_s=0\}$, conditional on $\{N_{0+}=k\}$ \tredd{describes the situation where all $k$ clients have left the system \textit{before} time $s$. Hence, this event} amounts to a Poisson process with intensity $\mu$ making {\it at least} $k$ jumps in the interval $[0,s]$. Likewise, $\{N_s=k-\ell\}$ (with $\ell=0,\ldots,k-1$), conditional on $\{N_{0+}=k\}$, \tredd{amounts to} this Poisson process making {\it precisely} $\ell$ jumps in $[0,s]$, \tredd{which is precisely the number of clients that have left in this interval}. We thus find
\begin{align}\label{fg}f_k(t) &=  \int_0^t \sum_{\ell = k}^\infty e^{-\mu s}\frac{(\mu s)^\ell}{\ell !}{\rm d}s,\hspace{1cm}
g_k(t)= \int_0^t\sum_{\ell=0}^{k-1} (k-\ell-1) e^{-\mu s}\frac{(\mu s)^\ell}{\ell !}{\rm d}s.
\end{align}

The above expressions for $f_k(t)$ and $g_k(t)$ can be considerably simplified. To this end, we first write the integrals featuring in \eqref{fg} in closed form. Define ${\rm E}(\ell,\mu)$ as an Erlang random variable with shape parameter $\ell$ and rate parameter $\mu.$ Then,
\begin{align}\int_0^t  e^{-\mu s}\frac{(\mu s)^\ell}{\ell !}\,{\rm d}s &=\frac{1}{\mu}\int_0^{\mu t} e^{-u}\frac{u^\ell}{\ell !}\,{\rm d}u =\mu^{-1} \cdot{\mathbb P}({\rm E}(\ell+1,1)\leqslant \mu t)=\frac{1}{\mu}\sum_{m=\ell+1}^\infty e^{-\mu t}\frac{(\mu t)^m}{m!}.\label{int}
\end{align}
Let us start by evaluating $f_k(t)$. By virtue of \eqref{int},
\begin{align*}f_k(t)&= \sum_{\ell = k}^\infty\frac{1}{\mu}\sum_{m=\ell+1}^\infty e^{-\mu t}\frac{(\mu t)^m}{m!}=
\frac{1}{\mu} \sum_{m=k+1}^\infty \sum_{\ell=k}^{m-1} e^{-\mu t}\frac{(\mu t)^m}{m!}=\frac{1}{\mu} \sum_{m=k+1}^\infty (m-k)e^{-\mu t}\frac{(\mu t)^m}{m!}.
\end{align*}
Hence,  with ${\rm Pois}(\mu)$ denoting a Poisson random variable with mean $\mu$, we   obtain
\begin{align*}f_k(t)&=\frac{1}{\mu} \sum_{m=k+1}^\infty  e^{-\mu t}\frac{(\mu t)^{m}}{(m-1)!}-\frac{k}{\mu}
\sum_{m=k+1}^\infty  e^{-\mu t}\frac{(\mu t)^{m}}{m!}\\
&={ t} \sum_{m=k}^\infty  e^{-\mu t}\frac{(\mu t)^{m}}{m!}-\frac{k}{\mu}
\sum_{m=k+1}^\infty  e^{-\mu t}\frac{(\mu t)^{m}}{m!} =  t\cdot {\mathbb P}({\rm Pois}(\mu t)\geqslant k)- \frac{k}{\mu} \cdot{\mathbb P}({\rm Pois}(\mu t)\geqslant k+1).
\end{align*}
The quantity $g_k(t)$ can be dealt with in a similar way, but the calculations are slightly more involved. Appealing to \eqref{int}, with $\ell':=k-\ell-1$,
\[g_k(t)=\sum_{\ell=0}^{k-1} (k-\ell-1)\left(\frac{1}{\mu}\sum_{m=\ell+1}^\infty e^{-\mu t}\frac{(\mu t)^m}{m!}\right)
=\frac{1}{\mu}\sum_{\ell'=0}^{k-1}  \ell'\sum_{m=k-\ell'}^\infty e^{-\mu t}\frac{(\mu t)^m}{m!}.\]
Swapping the order of the sums, the right-hand side of the previous display can be alternatively written as
\begin{align}\frac{1}{\mu}&\sum_{m=1}^k e^{-\mu t}\frac{(\mu t)^m}{m!}\sum_{\ell'=k-m}^{k-1}  \ell'+
\frac{1}{\mu}\sum_{m=k+1}^\infty e^{-\mu t}\frac{(\mu t)^m}{m!}\sum_{\ell'=0}^{k-1}  \ell' \nonumber \\
&=\frac{1}{\mu}\sum_{m=1}^k e^{-\mu t}\frac{(\mu t)^m}{m!} \left({{k}\choose{2}}-{{k-m}\choose{2}}\right)+
\frac{1}{\mu}\sum_{m=k+1}^\infty e^{-\mu t}\frac{(\mu t)^m}{m!} {{k}\choose{2}}\nonumber \\
&=\frac{1-e^{-\mu t}}{\mu}{{k}\choose{2}}-\frac{1}{\mu}\sum_{m=1}^k e^{-\mu t}\frac{(\mu t)^m}{m!} {{k-m}\choose{2}}. \label{G1}
\end{align}
We proceed by analyzing the second term in \eqref{G1}.
It is easily verified that
${{k-m}\choose{2}}$ can be written as $\tfrac{1}{2}k^2 -m(k-1) +\tfrac{1}{2} m(m-1) -\tfrac{1}{2}k.$
An application of this identity directly leads to, with the empty sum being defined as $0$,
\[\sum_{m=1}^k  \frac{(\mu t)^m}{m!} {{k-m}\choose{2}}=\frac{k(k-1)}{2}\sum_{m=1}^k  \frac{(\mu t)^m}{m!}
-(k-1)\sum_{m=1}^k  \frac{(\mu t)^m}{(m-1)!}+\frac{1}{2}
\sum_{m=2}^k  \frac{(\mu t)^m}{(m-2)!}.\]
This implies that 
\begin{align*}\frac{1}{\mu}\sum_{m=1}^k e^{-\mu t}\frac{(\mu t)^m}{m!} {{k-m}\choose{2}}&=\frac{k(k-1)}{2\mu}\cdot {\mathbb P}(1\leqslant {\rm Pois}(\mu t)\leqslant k)\:-\\
&\:\:\:(k-1)t \cdot {\mathbb P}( {\rm Pois}(\mu t)\leqslant k-1)+\frac{\mu t^2}{2} \cdot {\mathbb P}( {\rm Pois}(\mu t)\leqslant k-2).
\end{align*}
Upon combining the above,
\begin{align*}
g_k(t)&=\frac{k(k-1)}{2\mu}\cdot {\mathbb P}({\rm Pois}(\mu t)\geqslant 1)-\frac{k(k-1)}{2\mu}\cdot {\mathbb P}(1\leqslant {\rm Pois}(\mu t)\leqslant k)\,+\\&\hspace{14mm}\,(k-1)t \cdot {\mathbb P}( {\rm Pois}(\mu t)\leqslant k-1)
-\frac{\mu t^2}{2} \cdot {\mathbb P}( {\rm Pois}(\mu t)\leqslant k-2).
\end{align*}

The following proposition summarizes what we have found above. We use the notation $F_\mu(k):= {\mathbb P}({\rm Pois}(\mu)\leqslant k)$ for the distribution function of a Poisson random variable with mean $\mu.$
\begin{proposition} \label{P1} For $k=1,\ldots,i$ and $t\geqslant 0$,
\begin{align*}f_k(t) &= t\cdot (1-F_{\mu t}(k-1))- \frac{k}{\mu} \cdot(1-F_{\mu t}(k)),\\
g_k(t)&=(k-1)t \cdot F_{\mu t}( k-1)
-\frac{\mu t^2}{2} \cdot F_{\mu t}(k-2)+\frac{k(k-1)}{2\mu}\cdot (1-F_{\mu t}(k)).
\end{align*}
\end{proposition}

\tredd{Besides the quantities $f_k(t)$ and $g_k(t)$, in the dynamic-programming routine an important role is played by the transition probabilities: for $k=1,\ldots,i$ and $\ell=1,\ldots,k+1$,}
\[p_{k\ell}(t):={\mathbb P}(N_{t+}=\ell\,|\,N_{0+}=k).\]
These can be computed in a standard way; we state their explicit form for completeness. 
\begin{lemma} \label{L1}For $k=1,\ldots,i$ and $\ell=2,\ldots,k+1$ and $t\geqslant 0$,
\begin{align*}p_{k1}(t)&= \sum_{m = k}^\infty e^{-\mu t}\frac{(\mu t)^m}{m !},\:\:\:\:
p_{k\ell}(t)=e^{-\mu t}\frac{(\mu t)^{k-\ell+1}}{(k-\ell+1) !}.
\end{align*}
\end{lemma}

We can now apply dynamic programming to find the optimally scheduled arrival time of the next client. Recall that $C_i(k)$,  with $i=1,\ldots,n$ and $k=1,\ldots,i$, corresponds to the cost incurred from the arrival of the $i$-th client, given there are $k$ clients in the system immediately after the arrival of this $i$-th client. {\textcolor{black}{This cost can be split up into two parts: the \textit{immediate cost}, induced by the decision of scheduling client $i+1$ precisely $t$ time units later, and the \textit{remaining cost} after this scheduled arrival. Note that client $n$ is the last client to arrive, after which no decision regarding a next arrival is to be made. Since at the arrival of this client, the next idling moment occurs when the service of all clients has been completed, the cost $C_n(k)$ only consists of the remaining waiting times of the clients in the system.}}
\begin{theorem} \label{t1}  Let $f_k(t)$, $g_k(t)$ and $p_{k\ell}(t)$ be given by Proposition $\ref{P1}$ and  Lemma $\ref{L1}$. 
We can determine the $C_i(k)$  recursively: for $i=1,\ldots,n-1$ and $k=1,\ldots,i$,
\[C_{i}(k) =\inf_{t\geqslant 0}\left(\omega\,f_k(t)+(1-\omega)\,g_k(t)+ \sum_{\ell=1}^{k+1} p_{k\ell}(t)\, C_{i+1}(\ell)\right),\]
whereas, for $k=1,\ldots,n$,
\[C_{n}(k) =(1-\omega)\,g_k(\infty)=(1-\omega)\,\frac{k(k-1)}{2\mu}.\]
\end{theorem}

\begin{remark}\label{rem:g_h}{\em 
\tredd{In Theorem \ref{t1},} we \tredd{kept track of} the contributions of the waiting times to the slots between subsequent arrivals. More precisely, in the quantity $g_k(t)$ we keep track of the contributions to the expected waiting time by all \tredd{clients} in the system, between two subsequent arrivals (with an interarrival time $t$), given there are $k$ clients present at the beginning of the slot (i.e., directly after the arrival of a client). There is a convenient alternative, though: we can work instead with $h_k$, denoting the expected waiting time of a client if the number of clients immediately after her arrival is $k$  (i.e., also the waiting time {\it outside} the interarrival time under consideration). \tredd{Some thought reveals that this way one also gathers the sum of the clients' mean  waiting times. Indeed, the quantities $g_k(t)$ represent the aggregate waiting time (of all clients present, that is) in the interval of length $t$, whereas $h_k$ is the waiting time of the single client in the rest of the schedule, in both cases conditional on $k$ clients being present at the beginning of the slot. Informally one could say that working with the $g_k(t)$\,s one first adds up the contributions of the clients and then one aggregates over time, whereas working with the $h_k$\,s one does the opposite.} 

Working with the $h_k$\,s has an important advantage, namely its simplicity: we have 
\begin{equation}
\label{gster}
h_k =\frac{k-1}{\mu},\end{equation}
\tredd{independently} of $t$. It leads to the following alternative dynamic programming scheme. Let $f_k(t)$ and $p_{k\ell}(t)$ be given by Proposition $\ref{P1}$ and Lemma $\ref{L1}$, and let $h_k$ be given by \eqref{gster}. 
Then we can determine, for $i=1,\ldots,n-1$ and $k=1,\ldots,i$, the \tredd{$\tilde{C}_i(k)$} recursively by iterating
\[\tredd{\tilde{C}_{i}(k) =\inf_{t\geqslant 0}\left(\omega\,f_k(t)+(1-\omega)\,h_k + \sum_{\ell=1}^{k+1} p_{k\ell}(t)\, \tilde{C}_{i+1}(\ell)\right),}\]
whereas, for $k=1,\ldots,n$,
\[\tredd{\tilde{C}_{n}(k) =(1-\omega)\,h_k =(1-\omega)\,\frac{k-1}{\mu}.}\]
The resulting value of the objective function (i.e., \tredd{$\tilde{C}_1(1)$}) and the corresponding schedule coincide with their counterparts that are produced by the algorithm of Theorem \ref{t1}.}$\hfill\Diamond$
\end{remark}

{\small
\begin{table}
 \centering
\begin{tabular}{>{}l<{}|*{14}{c}}
\multicolumn{1}{l}{$i$} &&&&&&&&\\\cline{1-1}
1 & 0.88 & & & & & & & & & & & & & \\
2 & 0.88 & 1.94 & & & & & & & & & & & & \\
3 & 0.88 & 1.94 & 2.99 & & & & & & & & & & & \\
4 & 0.88 & 1.94 & 2.99 & 4.03 & & & & & & & & & & \\
5 & 0.88 & 1.94 & 2.99 & 4.03 & 5.06 & & & & & & & & & \\
6 & 0.88 & 1.94 & 2.99 & 4.03 & 5.06 & 6.09 & & & & & & & & \\
7 & 0.88 & 1.94 & 2.99 & 4.03 & 5.06 & 6.09 & 7.11 & & & & & & & \\
8 & 0.88 & 1.94 & 2.99 & 4.03 & 5.06 & 6.09 & 7.11 & 8.14 & & & & & & \\
9 & 0.88 & 1.94 & 2.99 & 4.03 & 5.06 & 6.09 & 7.11 & 8.14 & 9.16 & & & & & \\
10 & 0.88 & 1.94 & 2.99 & 4.03 & 5.06 & 6.09 & 7.11 & 8.14 & 9.16 & 10.18 & & & & \\
11 & 0.88 & 1.94 & 2.99 & 4.03 & 5.06 & 6.09 & 7.11 & 8.14 & 9.16 & 10.18 & 11.19 & & & \\
12 & 0.88 & 1.94 & 2.99 & 4.03 & 5.06 & 6.09 & 7.11 & 8.13 & 9.15 & 10.17 & 11.19 & 12.21 & & \\
13 & 0.86 & 1.91 & 2.96 & 3.99 & 5.02 & 6.04 & 7.07 & 8.09 & 9.11 & 10.12 & 11.14 & 12.15 & 13.17 & \\
14 & 0.69 & 1.68 & 2.67 & 3.67 & 4.67 & 5.67 & 6.67 & 7.67 & 8.67 & 9.67 & 10.67 & 11.67 & 12.67 & 13.67 \\\hline
\multicolumn{1}{l}{} &1&2&3&4&5&6&7&8&9&10&11&12&13&14\\\cline{2-15}
\multicolumn{1}{l}{} &\multicolumn{14}{c}{$k$}
\end{tabular}\caption{Optimal interarrival times $\tau_i(k)$ for Example~\ref{E0}. \label{TAB1}}
\end{table}}

\begin{example}\label{E0}{\em 
We present an instance in which we apply Theorem \ref{t1}. Denote by
\[\tau_i(k) := \arg\inf_{t\geqslant 0}\left(\omega\,f_k(t)+(1-\omega)\,g_k(t)+ \sum_{\ell=1}^{k+1} p_{k\ell}(t)\, C_{i+1}(\ell)\right)\]
the optimal time between the $i$-th and $(i+1)$-st arrival, given that immediately after the $i$-th arrival there are $k=1,\ldots,i$ clients present. 
Consider the case of $\omega = \frac{1}{2}$ and $n=15.$ We assume $\mu=1$, but values corresponding to any other value of $\mu$ can be found by \tredd{multiplying all interarrival times} by $1/\mu.$ The optimal interarrival times $\tau_i(k)$ corresponding to this schedule are given in Table \ref{TAB1}.
}$\hfill\Diamond$
\end{example}

\begin{example}\label{E1}{\em In this example, we assess the gain achieved by scheduling dynamically rather than in advance. We let $K_{\rm pre}(n,\omega)$ be the value of the objective function (for a given $\omega$) when we precalculate the schedule, i.e., optimizing \eqref{OBJ} at time $0$ over all arrival times $t_1=0,t_2,\ldots,t_n$. Also, $K_{\rm dyn}(n,\omega)$ is the value of the objective function when at every arrival we determine the optimal arrival epoch for the next client, taking into account the number of clients present. In addition, we define $r(n,\omega)$ as the ratio of $K_{\rm dyn}(n,\omega)$ and $K_{\rm pre}(n,\omega)$.
The results, shown in Table~\ref{TAB2}, indicate that the gain of scheduling dynamically can be substantial. In particular, this holds when the weight $\omega$ is relatively high and/or the number of clients $n$ gets large. \hfill$\Diamond$
}
\end{example}

\begin{remark}\label{RRR}
{\em We observed that the gain of adapting the schedule becomes more pronounced as $\omega$ increases. This phenomenon is essentially due to the intrinsic asymmetry between idle times and waiting times. 
    In the regime that $\omega$ is relatively large, the schedule will be such that interarrival times are relatively short, so as to avoid idle time. This means that there are systematically relatively many clients in the system, where for each number of clients there is a specific ideal (i.e., cost-minimizing) time \tredd{until} the next arrival. In other words, there is a substantial gain due to the fact that current information is taken into account when determining the next arrival time. If $\omega$ is relatively small, on the contrary, interarrival times are relatively long to avoid waiting times. In the extreme case of a very small $\omega$, effectively every new \tredd{client} will find the system empty. This means that there is no gain in letting the time until the next arrival depend on the state of the system, as the state of the system upon client arrival is virtually always the same. \hfill$\Diamond$ }
\end{remark}

{\small
\begin{table}
\centering
\begin{tabular}{|c|c|ccccccccc|}
\hline
$n$   & $\omega$        & 0.1  & 0.2  & 0.3  & 0.4  & 0.5  & 0.6  & 0.7  & 0.8  & 0.9 \\
\hline
5     & $K_{\rm dyn}(n,\omega)$ & 0.94 & 1.36 & 1.58 & 1.67 & 1.65 & 1.54 & 1.34 & 1.04 & 0.61 \\
      & $K_{\rm pre}(n,\omega)$ & 0.98 & 1.46 & 1.74 & 1.87 & 1.88 & 1.78 & 1.56 & 1.21 & 0.71 \\
      & $r(n,\omega)$          & 0.96 & 0.93 & 0.91 & 0.89 & 0.88 & 0.87 & 0.86 & 0.86 & 0.86 \\
\hline
10    & $K_{\rm dyn}(n,\omega)$ & 2.13 & 3.09 & 3.62 & 3.85 & 3.85 & 3.64 & 3.21 & 2.55 & 1.60 \\
      & $K_{\rm pre}(n,\omega)$ & 2.25 & 3.39 & 4.12 & 4.54 & 4.69 & 4.58 & 4.19 & 3.44 & 2.21 \\
      & $r(n,\omega)$          & 0.95 & 0.91 & 0.88 & 0.85 & 0.82 & 0.79 & 0.77 & 0.74 & 0.72 \\
\hline
15    & $K_{\rm dyn}(n,\omega)$ & 3.32 & 4.83 & 5.66 & 6.04 & 6.05 & 5.73 & 5.08 & 4.07 & 2.57 \\
      & $K_{\rm pre}(n,\omega)$ & 3.51 & 5.33 & 6.51 & 7.23 & 7.55 & 7.47 & 6.94 & 5.85 & 3.92 \\
      & $r(n,\omega)$          & 0.95 & 0.91 & 0.87 & 0.83 & 0.80 & 0.77 & 0.73 & 0.70 & 0.66 \\
\hline
20    & $K_{\rm dyn}(n,\omega)$ & 4.51 & 6.56 & 7.70 & 8.22 & 8.25 & 7.83 & 6.96 & 5.58 & 3.54 \\
      & $K_{\rm pre}(n,\omega)$ & 4.78 & 7.27 & 8.90 & 9.93 & 10.41& 10.36& 9.72 & 8.32 & 5.73 \\
      & $r(n,\omega)$          & 0.95 & 0.90 & 0.87 & 0.83 & 0.79 & 0.76 & 0.72 & 0.67 & 0.62 \\
\hline
25    & $K_{\rm dyn}(n,\omega)$ & 5.70 & 8.29 & 9.74 & 10.40& 10.45& 9.92 & 8.83 & 7.09 & 4.51 \\
      & $K_{\rm pre}(n,\omega)$ & 6.04 & 9.21 & 11.30& 12.62& 13.28& 13.27& 12.52& 10.82& 7.60 \\
      & $r(n,\omega)$          & 0.94 & 0.90 & 0.86 & 0.82 & 0.79 & 0.75 & 0.71 & 0.66 & 0.59 \\
\hline
30    & $K_{\rm dyn}(n,\omega)$ & 6.89 & 10.03& 11.77& 12.59& 12.65& 12.02& 10.70& 8.61 & 5.48 \\
      & $K_{\rm pre}(n,\omega)$ & 7.30 & 11.14& 13.69& 15.32& 16.14& 16.18& 15.32& 13.33& 9.50 \\
      & $r(n,\omega)$          & 0.94 & 0.90 & 0.86 & 0.82 & 0.78 & 0.74 & 0.70 & 0.65 & 0.58 \\
\hline
\end{tabular}
\caption{Cost of dynamic and precalculated schedule for Example~\ref{E1}. \label{TAB2}}
\end{table}}

\tred{We conclude this section discussing the concept of {\it stationary schedules.} Table \ref{TAB1} reveals that the $\tau_i(k)$ depends on $k$ (i.e., the number of clients present), but hardly on $i$ (i.e., the index of the client who just arrived). Indeed, this suggests that there is an `approximate stationary policy', that provides $\tau(k)$, i.e, the interarrival time  when $k$ clients are present in settings where the total number of clients $n$ is large. Such a policy can be evaluated as follows. }

\tred{
Suppose that, when observing $k$ clients after a new client has joined, we consistently schedule the next client after $x_k$ time units. Then the number of clients immediately after client arrivals forms a (discrete-time) Markov chain, with the transition probability of going from $k$ to $\ell$ clients given by $p_{k\ell}(x_k)$ (where $p_{k\ell}(x)$ has been defined above). The equilibrium distribution of the Markov chain is given by ${\boldsymbol \pi}=(\pi_1,\pi_2,\ldots)$. Observe that ${\boldsymbol \pi}$ depends on all the interarrival times $x_k$ (due to the fact that the Markov chain depends on all the $x_k$), so that we prefer to write ${\boldsymbol \pi}({\boldsymbol x})$, with ${\boldsymbol x}=(x_1,x_2,\ldots)$. Then the minimal long-term per-client cost is given by
\[C = \inf_{\boldsymbol x} \left(\sum_{k=1}^\infty \pi_{k}({\boldsymbol x}) \big(
\omega\,f_k(x_k)+(1-\omega)\,g_k(x_k)
\big)\right),\]
with $f_k(x)$ and $g_k(x)$ as defined above.
Then $\tau(k)$ is the $k$-th component of the optimizing ${\boldsymbol x}.$
To numerically evaluate the long-term per-client cost $C$, there is the practical (but not complicated) issue that the above sum should be truncated at a suitably chosen level (say $K$). }

\tred{It is noted that the $\tau(k)$ that are found are of great practical interest, as these can serve as the initial values for the $\tau_i(k)$ in the numerical evaluation of our (non-stationary) dynamic-programming based schedule, thus accelerating computation times considerably. The search for the optimizing ${\boldsymbol x}$ is of low computational cost, as this concerns just the minimization of a function with a $K$-dimensional argument, which can be done relying on standard software. It should be kept in mind that in every step we have to compute, for the current value of ${\boldsymbol x}=(x_1,\ldots,x_K)$, the equilibrium distribution ${\boldsymbol \pi}({\boldsymbol x})$ corresponding to the transition probabilities $p_{k\ell}(x_k)$, but this is a matter of using a standard routine to solve a system of linear equations. In Table \ref{TAB2a} we present the optimal stationary dynamic schedules for various values of the weight $\omega$. Notice the striking agreement between the $\tau_i(k)$ values presented in Table \ref{TAB1} on one hand, and the $\tau(k)$ values in the column corresponding to $\omega=0.5$ in Table \ref{TAB2a} on the other hand.}

{\small
\begin{table}
\centering
\begin{tabular}{|c|ccccccccc|}
\hline
 $\omega$       & 0.1  & 0.2  & 0.3  & 0.4  & 0.5  & 0.6  & 0.7  & 0.8  & 0.9 \\
\hline $k=1$ & 2.38& 1.73& 1.36& 1.09& 0.88& 0.70& 0.53& 0.38& 0.22\\
 $k=2$ & 3.98& 3.15& 2.64& 2.26& 1.94& 1.66& 1.39& 1.10& 0.77\\
  $k=3$ & 5.42& 4.45& 3.85& 3.39& 2.99& 2.63& 2.28& 1.90& 1.44\\
  $k=4$ & 6.79& 5.71& 5.02& 4.49& 4.03& 3.60& 3.18& 2.72& 2.15\\
  $k=5$ & 8.11& 6.93& 6.17& 5.58& 5.06& 4.58& 4.10& 3.57& 2.90\\
  $k=6$ & 9.40& 8.12& 7.30& 6.65& 6.09& 5.56& 5.02& 4.43& 3.66\\
\hline
\end{tabular}
\caption{Optimal interarrival times $\tau(k)$ corresponding to the stationary schedule. \label{TAB2a}}
\end{table}}

\tred{Stationary schedules have been studied in the static case as well; see e.g.\ \cite{KUIP}, and also the related heavy-traffic analysis of \cite{KMM}.}

\section{Heterogeneous exponential case} 
\label{hetexp} In this section, we consider the case that the service requirements are heterogeneous, but still exponentially distributed. As will become clear in the next section, some of the elements appearing in the present section  help in the analysis of the case of service times with a general SCV, as dealt with in the next section. 
\tredd{In particular, our approach deals with evaluating convolutions and recursions, which is thoroughly used in the general case in the next section. The main result of this section is Theorem~\ref{t2}, which is the equivalent of Theorem~\ref{t1} for heterogeneous exponential service times.}

The mean of $B_i$, denoting the service requirement of client $i\in\{1,\ldots,n\}$, is $\mu_i^{-1}\in(0,\infty)$. For ease we assume that all $\mu_i$ are distinct; later in the section, we comment on the case that some of the $\mu_i$ coincide. 
In the sequel, we will intensively use, for $k=1,2,\ldots$, $\ell=0,1,\ldots$ and $s\geqslant 0$, the following notation for the density of the sum of independent exponentially distributed random variables:
\[\varphi_{k\ell}(s):= \frac{{\rm d}}{{\rm d}s}{\mathbb P}\left(\sum_{j=k}^{k+\ell} {\rm E}_j \leqslant s\right).\]
Here ${\rm E}_j$ denotes an exponentially distributed random variable with mean $\mu_j^{-1}$. The following lemma shows that this density can be written as a mixture of exponential terms, with the $\mu_j$ in the exponent. It should be noted that the corresponding weights are not necessarily positive. It is in principle possible to derive a closed-form expression for $\varphi_{k\ell}(s)$; this result is included in \tredd{Appendix \ref{app:proofs}}.

\begin{lemma} \label{L2} For $k=1,2,\ldots$ and $\ell=0,1,\ldots$ and $s\geqslant 0$, there are constants $c_{k\ell j}\in {\mathbb R}$, with $j=k,\ldots,k+\ell$, such that 
\[\varphi_{k\ell}(s) = \sum_{j=k}^{k+\ell} c_{k\ell j} \,e^{-\mu_j s}.\]
The coefficients $c_{k\ell j}$ are given recursively through $c_{k0k}=\mu_k$ and
\[ c_{k,\ell+1, j}=c_{k\ell j} \frac{\mu_{k+\ell+1}}{\mu_{k+\ell+1}-\mu_j }\:\:\:\mbox{for $j=k,\ldots,k+\ell$},\:\:\:\:\:\:c_{k,\ell+1, k+\ell+1}= \sum_{j=k}^{k+\ell} c_{k\ell j} \frac{\mu_{k+\ell+1}}{\mu_j-\mu_{k+\ell+1} }.\]
\end{lemma}

{\it Proof}: We prove the statement by induction. It is clear that $c_{k0k}=\mu_k$, as claimed, due to $\varphi_{k0}(s)=\mu_k \,e^{-\mu_k s}$ for $s\geqslant 0$. We proceed by verifying the induction step, by checking whether $\varphi_{k,\ell+1}(s)$ has the right form, given that $\varphi_{k\ell}(s)$ has. To this end, note that using the usual convolution representation, 
\begin{equation}\label{conv_rep}\varphi_{k,\ell+1}(s)= \int_0^s \varphi_{k\ell}(u) \,\mu_{k+\ell+1}e^{-\mu_{k+\ell+1}(s-u)}{\rm d}u,\end{equation}
which by the induction hypothesis can be written as
\begin{equation}\label{erl} \int_0^s  \sum_{j=k}^{k+\ell} c_{k\ell j} \,e^{-\mu_j u} \,\mu_{k+\ell+1}e^{-\mu_{k+\ell+1}(s-u)}{\rm d}u=
\sum_{j=k}^{k+\ell} c_{k\ell j} \frac{\mu_{k+\ell+1}}{\mu_j-\mu_{k+\ell+1} }(e^{-\mu_{k+\ell+1} s} - e^{-\mu_{j} s}).\end{equation}
This expression reveals the stated recursion. $\hfill\Box$

\begin{remark}{\em In the above setup we assumed that all $\mu_i$ are distinct, but the calculations can be adapted to the case that some $\mu_i$ coincide. Upon inspecting the proof of Lemma \ref{L2}, it is seen that the terms in the expression for $\varphi_{k\ell}(s)$ will not be just exponentials, but products of polynomials and exponentials. More concretely, if  in the right-hand side of \eqref{erl}, for some $j$ we have that $\mu_{k+\ell+1} = \mu_j$, we obtain that the corresponding summand should read 
\[c_{k\ell j} \mu_j s e^{-\mu_j s}.\]
It can be checked that if $m$ of the $\mu_j$ are equal, this leads to a term proportional to $s^{m-1} e^{-\mu_j s}$; cf.\ the density of the Erlang distribution. The computations below are performed for the case that
$\varphi_{k\ell}(s)$ is the sum of exponentials, but it can easily be adapted to the case that some terms are products of polynomials and exponentials.
$\hfill\Diamond$
}\end{remark}

The following interesting result is proved in \tredd{Appendix \ref{app:proofs}} \tredd{and will be useful to simplify several expressions later in this section.}

\begin{lemma} \label{lem:sum_to_1} For $k=1,\dots,n$ and $\ell = 0,\dots,n-k$,
\[
\sum_{j=k}^{k+\ell}\frac{c_{k\ell j}}{\mu_j} = 1.
\]
\end{lemma}

In this case with heterogeneous exponentially distributed service times, the setup of its homogeneous counterpart essentially carries over, but now  we have to work with the \tredd{terms}
\[f_{ki}(t):= \int_0^t {\mathbb P}_i\big(N_s=0\,|\, N_{0+}=k\big) {\rm d}s,\:\:\:g_{ki}(t):=\int_0^t\sum_{\ell=0}^{k-1} (k-\ell-1) {\mathbb P}_i\big(N_s=k-\ell\,|\, N_{0+}=k\big) {\rm d}s,\]
where the subscript $i$ indicates that the $i$-th client entered at time $0$. Again, to make the notation lighter, we shift in these computations time such that the arrival of the $i$-th client corresponds to time $0$. 

We start our computations with the evaluation of $f_{ki}(t)$, by first considering the integrand. It is readily checked that, by Lemma \ref{L2},
\begin{align*}
{\mathbb P}_i\big(N_s=0\,|\, N_{0+}=k\big) = {\mathbb P}\left(\sum_{j=i-k+1}^i {\rm E}_j \leqslant s\right)
= \sum_{j=i-k+1}^i c_{i-k+1,k-1,j}\,\int_{0}^{s}e^{-\mu_j u}{\rm d}u.
\end{align*}
Observe that the event of interest corresponds to clients $i-k+1$ up to (and including) $i$ being served before time $s$. By Lemma \ref{lem:sum_to_1}, we obtain
\[f_{ki}(t)
=\sum_{j=i-k+1}^i c_{i-k+1,k-1,j}\int_0^t\psi_j(s){\rm d}s=
t - \sum_{j=i-k+1}^i \frac{c_{i-k+1,k-1,j}}{\mu_j}\psi_j(t),\]
with
\[\psi_j(t):= \int_{0}^{t}e^{-\mu_j s} {\rm d}s = \frac{1-e^{-\mu_j t}}{\mu_j}.\]
As a sanity check, observe that indeed $f_{ki}(0)=0$, as desired: the mean amount of idle time over a period of length $0$ should equal $0$. Likewise, $f_{ki}(\infty)=\infty$.

We proceed by analyzing $g_{ki}(t)$. Using \eqref{conv_rep}, we observe that, for $\ell=0,\ldots,k-1$,
\begin{align*}
{\mathbb P}_i\big(N_s=k-\ell\,|\, N_{0+}=k\big) &=
\int_{0}^{s}\mathbb{P}\Bigg(\sum_{j=i-k+1}^{i-k+\ell}E_j \in \mathrm{d}u \Bigg)\mathbb{P}\left(E_{i-k+\ell+1} > s - u\right)\mathrm{d}u \\
&= \int_{0}^{s}\varphi_{i-k+1,\ell-1}(u)e^{-\mu_{i-k+\ell+1}(s-u)}\mathrm{d}u
= \frac{\varphi_{i-k+1,\ell}(s)}{\mu_{i-k+\ell+1}}.
\end{align*}
Now the event of interest corresponds to clients $i-k+1$ up to (and including) $i-k+\ell$ being served before time $s$, but the service of client $i-k+\ell+1$ being completed only after time $s$. It follows that
\[
g_{ki}(t)
= \sum_{\ell=0}^{k-1}(k-\ell-1)\sum_{j=i-k+1}^{i-k+\ell+1}\frac{c_{i-k+1,\ell,j}}{\mu_{i-k+\ell+1}}\psi_{j}(t).
\]

Now that we have expressions for the mean idle times and service times, we compute the corresponding transition probabilities.
Let these transition probabilities be defined as $p_{k\ell, i}(t):={\mathbb P}_i(N_{t+}=\ell\,|\,N_{0+}=k).$ They can be calculated using the machinery that we developed above: for $\ell=2,\ldots,k+1$,
\[
p_{k1,i}(t) = 1 - \sum_{\ell=2}^{k+1}p_{k\ell,i}(t),\quad
p_{k\ell,i}(t) = \frac{\varphi_{i-k+1,k-\ell+1}(t)}{\mu_{i-\ell+2}}.
\]

As for the homogeneous case, we can set up a dynamic program to identify the optimal strategy. 

\begin{theorem} \label{t2}  Let $f_{ki}(t)$, $g_{ki}(t)$ and $p_{k\ell, i}(t)$ be as defined above. 
We can determine the $C_i(k)$  recursively: for $i=1,\ldots,n-1$ and $k=1,\ldots,i$,
\[C_{i}(k) =\inf_{t\geqslant 0}\left(\omega\,f_{ki}(t)+(1-\omega)\,g_{ki}(t)+ \sum_{\ell=1}^{k+1} p_{k\ell,i}(t)\, C_{i+1}(\ell)\right),\]
whereas, for $k=1,\ldots,n$,
\[C_{n}(k) =(1-\omega)\,g_{kn}(\infty)=(1-\omega)\sum_{\ell=0}^{k-1} (k-\ell-1)\frac{1}{\mu_{n-k+\ell+1}}.\]\end{theorem}
As in the homogeneous case, in the above dynamic programming recursion, we can work with $h_{ki}$ instead of $g_{ki}(t)$, with now
\[h_{ki} :=\sum_{\ell=1}^{k-1} \frac{1}{\mu_{i-k+\ell}}.\]
With the above dynamic programming algorithm at our disposal, we can assess which order of the clients minimizes the objective function $C_1(1)$. 
\tredd{
Since this is not the main focus point of the paper, we refer the interested reader to Appendix~\ref{app:hetexp}, where the numerical examples for this section can be found.
}


\section{Service times with general SCV}\label{phase}
Having dealt with the exponential case, in the present section we extend our analysis to the case of the service times having a general SCV.
This we do relying on the concept of 
phase-type distributions. The main idea behind our approach is that we fit our service-time distribution, which is characterized by its mean and SCV, to specific types of phase-type distributions. For these, we point out how to produce dynamic schedules using dynamic programming. 
In this section we consider the homogeneous phase-type case (i.e., the clients' service times are i.i.d.\ random variables with a given mean and SCV); at the expense of a substantial amount of additional calculations, this can in principle be extended to the corresponding heterogeneous case. 

We define \cite[Section III.4]{AS} phase-type distributions in the following way. Let ${\bs \alpha}\in{\mathbb R}^d$ be a probability vector, i.e., its entries are non-negative and sum to $1$. In addition, we have the transition rate matrix
\[Q = \left(\begin{array}{cc}T&{\bs t}\\
{\bs 0}^{\top}&0\end{array}\right),\]
with ${\bs 0}$  denoting a $d$-dimensional all zeroes (column-)vector, ${\bs t}:=-T{\bs 1}$, and ${\bs 1}$ a $d$-dimensional all ones (column-)vector. The pair $({\bs \alpha},T)$ represents a phase-type distributed random variable $B$: the initial phase is sampled from the distribution ${\bs \alpha}$, then the state evolves according to a continuous-time Markov chain with rate matrix $Q$, until the absorbing state $d+1$ is reached. The value of the phase-type distributed random variable $B$ records the time it takes for the process we just described to reach state $d+1$. The exponentially distributed times spent in the states $1,\ldots,d$ are often referred to as {\it phases}.

In the context of appointment scheduling with phase-type service requirements, an evident major complication is that the phase is a {\it non-observable} quantity. Indeed, when scheduling the arrival epoch of the $(i+1)$-st client, which happens at the moment client $i$ arrives, the information one has is the number of clients in the system and the elapsed service time of the client in service (in addition to the number of clients still to be served). Importantly, one does {\it not} know the phase that the client in service is in.  This can be resolved, however, by working with the distribution of the phase, conditional on the elapsed service time. It is this crucial idea that is applied extensively in the procedure outlined below.

Two special subclasses of phase-type distributions are of particular interest: the weighted Erlang distribution and the hyperexponential distribution.
More concretely, as motivated in great detail in e.g.\ \cite{KUIP,tijms}, we propose to approximate the service-time
distribution by a weighted Erlang distribution in case the SCV is smaller than 1, whereas the hyperexponential distribution can be used if the SCV is larger than 1. \tred{In this context we refer to the {\it method of moments}, frequently used in statistics; see also \cite{ANO} for statistical procedures to estimate phase-type distributions from data.}

In the next subsection, we detail the procedure to identify the parameters of the weighted Erlang distribution and the hyperexponential distribution. 
Importantly, the error due to replacing a non-phase-type distribution with its phase-type counterpart (with the same mean and SCV) is negligible, as was shown in e.g.\ \cite[pp.\ 110-111]{AThesis}. \tred{We provide a detailed account of this aspect in Section \ref{sec:discu}.}

\subsection{Fitting phase-type distributions}\label{sec:phase_type_fit}

\noindent As mentioned above, we aim at identifying, for any given distribution of the service time $B$,  a phase-type distribution that matches the first and second moment, or, equivalently, the mean and the SCV. Here, the \textit{coefficient of variation} ${\mathbb C}{\rm V}(X)$ of a random variable $X$ is defined as the ratio of the standard deviation to the mean:
\[
{\mathbb C}{\rm V}(X) := \frac{\sqrt{\mathbb{E}(X - \mathbb{E}X)^2}}{\mathbb{E}X}.
\]
The \textit{squared coefficient of variation} ${\mathbb S}(X)$ \tredd{or SCV} is defined as the square of ${\mathbb C}{\rm V}(X)$:
\[
{\mathbb S}(X) := \frac{{\mathbb V}{\rm ar}(X)}{\mathbb{E}[X]^2} = \frac{\mathbb{E}X^2}{\mathbb{E}[X]^2} - 1.
\]

\noindent As in \cite{KUIP}, if  ${\mathbb S}(B)$ is smaller than 1, then we match a mixture of two Erlang distributions (also referred to as a weighted Erlang distribution) with $K$ and $K+1$ phases  (for a suitably chosen value of $K\in{\mathbb N}$).
If, on the contrary, ${\mathbb S}(B)$ is larger than 1, then we
use a mixture of two exponential distributions (also referred to as a hyperexponential distribution). When ${\mathbb S}(B)$  equals $1$, both cases reduce to the homogeneous exponential case considered in Section \ref{homexp}. 
\begin{itemize}
\item[$\circ$]
In the first case, i.e., when ${\mathbb S}(B) \leqslant 1$, our fitted  distribution can be represented as follows: with the obvious independence assumptions, for some $K\in {\mathbb N}$, $\mu>0$, and $p\in [0,1]$,
\[
B \sim \text{E}(K,\mu) \mathbbm{1}_{\{U < p\}} + \text{E}(K+1,\mu) \mathbbm{1}_{\{U > p\}},
\]
where $U \sim \text{Unif}[0,1]$. In words: with probability $p$ we have an Erlang distribution with $K$ phases and with probability $1-p$ an Erlang distribution with $K+1$ phases. As
\[
\mathbb{E}B = p\frac{K}{\mu} + (1-p)\frac{K+1}{\mu},\quad
\mathbb{E}B^2 = p\frac{K(K+1)}{\mu^2} + (1-p)\frac{(K+1)(K+2)}{\mu^2},
\]
we find $\mu = (K + 1 - p)/\mathbb{E}B$ and
\[
\frac{\mathbb{E}B^2}{\mathbb{E}[B]^2} = \frac{(K+1)(pK + (1-p)(K+2))}{(pK + (1-p)(K+1))^2}
= \frac{(K+1)(K + 2(1-p))}{(K+1-p)^2}.
\]
Hence,
\[
{\mathbb S}(B)
= \frac{(K+1)(K + 2(1-p) - (K+1) + 2p) - p^2}{(K+1-p)^2}
= \frac{K+1 - p^2}{(K+1-p)^2}.
\]
Now define the function $f(\cdot)$ (and its  derivative) through
\[
f(x) := \frac{K + 1 - x^2}{(K + 1 - x)^2},\quad
f'(x) = \frac{2(K + 1)(1 - x)}{(K + 1 - x)^3} > 0,\quad x \in [0,1].
\]
Due to $f(0) = 1/(K + 1)$ and $f(1) = 1/K$, we find that ${\mathbb S}(B)$ lies between these two values. Hence, given ${\mathbb S}(B) \leqslant 1$, we can set $K = \lfloor {1}/{{\mathbb S}(B)}\rfloor$. It remains to identify the probability $p$. Solving the equation
\[
({\mathbb S}(B) + 1)p^2 - 2(K+1){\mathbb S}(B)p + (K+1)^2{\mathbb S}(B) - (K+1) = 0
\]
yields
\[
p = \frac{(K+1) {\mathbb S}(B) \pm \sqrt{(K+1) (1 - K\cdot {\mathbb S}(B))}}{{\mathbb S}(B) + 1},
\]
the lowest of which (the solution with the minus sign, that is) lies in the interval $[0,1]$.
\item[$\circ$] In the second case, i.e., when ${\mathbb S}(B) > 1$, we fit a hyperexponential distribution:
\[
B \sim \text{Exp}(\mu_1)\mathbbm{1}_{\{U<p\}} + \text{Exp}(\mu_2)\mathbbm{1}_{\{U>p\}},
\]
where $\mu_1,\mu_2 > 0$, $\mu_1 > \mu_2$ (without loss of generality) and $U \sim \text{Unif}[0,1]$. We have
\[
\mathbb{E}B = p\frac{1}{\mu_1} + (1-p)\frac{1}{\mu_2},
\quad
\mathbb{E}B^2 = p\frac{2}{\mu_1^2} + (1-p)\frac{2}{\mu_2^2}.
\]
Note that there are now three parameters to pick, which have to satisfy two equations. To find a unique solution, as pointed out in \cite{KUIP,tijms}, we impose the additional condition of \textit{balanced means}, i.e.,
$\mu_1 = 2p\mu$ and $\mu_2 = 2(1-p)\mu$ for some $\mu > 0$. Then it follows that
\[
\mathbb{E}B = p\frac{1}{2p\mu} + (1-p)\frac{1}{2(1-p)\mu} = \frac{1}{\mu}.
\]
Hence, given the expectation $\mathbb{E}B$ of our service-time distribution, we find the parameters $\mu_1$ and $\mu_2$ by setting $\mu = 1/\mathbb{E}B$, and using the \textit{balanced means} condition. We thus find
\[
\frac{\mathbb{E}B^2}{\mathbb{E}[B]^2} = \mu^2 \left(p\frac{2}{4p^2 \mu^2} + (1 - p)\frac{2}{4(1-p)^2 \mu^2}\right)
= \frac{1}{2p} + \frac{1}{2(1-p)} = \frac{1}{2p(1-p)}.
\]
As a consequence, 
\begin{align*}
{\mathbb S}(B) =\frac{1}{2p(1-p)} - 1
&\iff
p^2 - p + \frac{1}{2({\mathbb S}(B) + 1)} = 0.
\end{align*}
Note that, as desired ${\mathbb S}(B) > 1$.
As $\mu_1 > \mu_2$, we obtain the unique solution
\[
p = \frac{1}{2}\left(1 + \sqrt{\frac{{\mathbb S}(B) - 1}{{\mathbb S}(B) + 1}}\,\right).
\]
\end{itemize}

\noindent In the following two subsections we develop dynamic programming algorithms for the weighted Erlang case and hyperexponential case, respectively.

\subsection{Weighted Erlang distribution}\label{WED}
In this case, the service time $B$ equals with probability $p\in[0,1]$ an Erlang-distributed random variable with $K$ exponentially distributed phases, each of them having mean $\mu^{-1}$, and with probability $1-p$ an Erlang-distributed random variable with $K+1$ exponentially distributed phases, again with mean $\mu^{-1}$. This means that the corresponding $T$-matrix, of dimension $(K+1)\times(K+1)$, is
\[T= \left(\begin{array}{ccccccc} -\mu&\mu&0&\cdots&0&0&0\\
0&-\mu&\mu&\cdots&0&0&0\\
0&0&-\mu&\cdots&0&0&0\\
\vdots&\vdots&\vdots&&\vdots&\vdots&\vdots\\
0&0&0&&-\mu&\mu&0\\
0&0&0&&0&-\mu&\mu(1-p)\\
0&0&0&&0&0&-\mu\\
\end{array}\right),\]
and $\alpha_1=1$\tredd{, i.e., ${\bs \alpha} = (1,0,\dots,0)^T$.}

Let $Z_t$ be the phase the process is in at time $t$. We argued above that the phase as such is not observable. For that reason, we wish to determine the distribution of the phase as a function of the (observable) elapsed service time. We compute for $t\geqslant 0$ and $z=1,\ldots,K+1$,
\[\gamma_z(t):={\mathbb P}(Z_t=z\,|\,B>t).\]
A routine calculation yields that $\gamma_z(t) =  \gamma^\circ_z(t)/\gamma^\circ(t)$, where
\begin{align*} \gamma^\circ(t)&:={\mathbb P}(B>t)
=\sum_{ z =1}^{K}e^{-\mu t}\frac{(\mu t)^{ z -1}}{( z -1)!}+
(1-p)e^{-\mu t}\frac{(\mu t)^K}{K!},\\
 \gamma^\circ_ z (t)&:={\mathbb P}(Z_t= z , B>t)= e^{-\mu t}\frac{(\mu t)^{ z -1}}{( z -1)!}\mathbbm{1}_{\{ z =1,\ldots,K\}}+(1-p)e^{-\mu t}\frac{(\mu t)^K}{K!}\mathbbm{1}_{\{ z =K+1\}}.
 \end{align*}
 It is easily verified that indeed $\gamma_1(0)=1$ and $\gamma_2(0)=\cdots=\gamma_{K+1}(0)=0$, as expected. Along the same lines, $\gamma_{K+1}(\infty)=1$ and $\gamma_1(\infty)=\cdots=\gamma_{K}(\infty)=0$.
 
 We proceed by successively evaluating the mean idle time, the mean waiting time, and the transition probabilities, this facilitating setting up the dynamic program. \tredd{We do this by first conditioning on the phase $z = 1,\dots,K+1$ of the client in service; for this we make use of the notation $\bar{f}_{kz}(t)$ and $\bar{g}_{kz}(t)$ (or equivalently, $\bar{h}_{kz}$); they can be seen as more general versions of the $f_k(t)$ and $g_k(t)$ (or $h_k$) of the homogeneous exponential case. Next, we focus on computing the mean idle and waiting time given the \textit{observable} information, i.e., the elapsed service time $u \geq 0$ besides the number of clients $k$ in the system. Since these terms depend on the newly introduced terms, we use the notation $\bar{f}^{\circ}_{k,u}(t)$ and $\bar{g}^{\circ}_{k,u}(t)$ (or $\bar{h}^{\circ}_{k,u}$) for the idle and waiting time given the state information.}

 \subsubsection*{$\rhd$ Mean idle time}
 We start by evaluating the mean idle time, given $Z_{0+}= z$. Using the same arguments as in the previous sections, this quantity equals
 \begin{align*}\bar f_{k z }(t)&:= \int_0^t {\mathbb P}\big(N_s=0\,|\, N_{0+}=k, Z_{0+}= z \big)\,{\rm d}s.\end{align*}
 Let ${\rm Bin}(z,p)$ denote a binomially distributed random variable with parameters $z\in{\mathbb N}$ and $p\in[0,1]$.
 If $ z  =1,\ldots,K$, to achieve $N_s=0$ we should have that at time $s$ the number of phases that have been served is
 $K- z +1$ (the `certain phases' of the client in service), plus 1 with probability $1-p$ (the `uncertain phase' of the client in service), plus $(k-1)K$ (the `certain phases' of all $k-1$ waiting clients), plus ${\rm Bin}(k-1,1-p)$ (the `uncertain phases' of all waiting clients). 
 Hence,
 \begin{align*}
 {\mathbb P}\big(N_s=0\,|\, N_{0+}=k ,Z_{0+}= z \big) &= \sum_{m=0}^k \mathbb{P}({\rm Bin}(k,1-p) = m)
 {\mathbb P}\big({\rm E}(kK- z +1 + m,\mu) \leqslant s\big).
 \end{align*}
As a consequence, with $f_k(t)$ as given by Proposition \ref{L1}, 
 \begin{equation}
 \label{C1}\bar f_{k z }(t) = \sum_{m=0}^k \mathbb{P}({\rm Bin}(k,1-p) = m) f_{kK- z +1 + m}(t).\end{equation}
The case that $ z =K+1$ has to be treated separately, as in this case, the number of phases corresponding to the client in service is $K+1$ with certainty. This means that for $z=K+1$, in order to achieve $N_s=0$ it is required  that at time $s$ the number of phases that have been served is~$1$ (the $(K+1)$-st phase of the client in service), plus $(k-1)K$ (the `certain phases' of all waiting clients), plus ${\rm Bin}(k-1,1-p)$ (the `uncertain phases' of all waiting clients). We thus arrive, for $z=K+1$, at
  \[{\mathbb P}\big(N_s=0\,|\, N_{0+}=k ,Z_{0+}= z \big) = \sum_{m=0}^{k-1} \mathbb{P}({\rm Bin}(k-1,1-p) = m)
 {\mathbb P}\big({\rm E}((k-1)K+1 + m,\mu) \leqslant s\big),\]
implying that
 \begin{equation}
 \label{C2}\bar f_{k z}(t) = \sum_{m=0}^{k-1} \mathbb{P}({\rm Bin}(k-1,1-p) = m) f_{(k-1)K+1 + m}(t).\end{equation}
Define by $\bar B_s$ the age of the service time (i.e., the elapsed service time) of the client who is in service at time $s\geqslant 0$. Recalling how we introduced $\gamma_z(u)$, it is evident that
\[\bar f^\circ_{k,u }(t):=  \int_0^t {\mathbb P}\big(N_s=0\,|\, N_{0+}=k ,\bar B_{0+}= u \big)\,{\rm d}s=
\sum_{z=1}^{K+1} \gamma_z(u) \bar f_{k z }(t),\]
where the $\bar f_{k z }(t)$ are given by \eqref{C1} for $z=1,\ldots,K$, and by \eqref{C2} for $z=K+1$. This completes our computation of the mean idle time given the observable quantities (viz.\ the number of clients present at the beginning of the slot and the elapsed service time of the client in service).

 \subsubsection*{$\rhd$ Mean waiting time}
Now that we have analyzed the mean idle times, we continue by computing the mean waiting times. \tredd{As pointed out in Remark \ref{rem:g_h}, it is more convenient to work with the term $\bar{h}^{\circ}_{k,u}$ representing the expected waiting time of the newly arrived client rather than $\bar{g}^{\circ}_{k,u}(t)$, i.e., the expected waiting time of all clients present between the client arrivals. For completeness, the evaluation of $\bar{g}^{\circ}_{k,u}(t)$ is provided in Appendix \ref{app:mean_idle_time}; we now focus on the computation of $\bar{h}^{\circ}_{k,u}$. We find}
 \[\bar h_{k,u}^\circ :=\sum_{z=1}^{K+1} \gamma_z(u)\bar h_{kz} ;\]
 with $\bar h_{1z}:=0$, and for $k=2,\ldots,i$,
 \[\bar h_{kz} :=\frac{(k-1)(K+1-p)+1-z}{\mu}\,1_{\{z=1,\ldots,K\}}+ \frac{(k-2)(K+1-p)+1}{\mu}\,1_{\{z=K+1\}}.\]

This completes our computation of the mean waiting time given the observable quantities.

 \subsubsection*{$\rhd$ Transition probabilities}
As a next step we evaluate the transition probabilities at arrival epochs, i.e., the distribution of $(N_{t+}, \bar B_{t+})$ conditional on $N_{0+}=k, \bar B_{0+}=u$. The analysis borrows elements from e.g.\ \cite{KEMP, W1, W2}. We start by separately analyzing the cases $N_{t+}=k+1$ and $N_{t+}=1$. \tredd{As these are the two most extreme scenarios (i.e., no client has left the system or all clients have left the system), we use the notation $P^\uparrow_{k,u}(t)$ and $P^\downarrow_{k,u}(t)$ for the corresponding transition probabilities.} First observe that for $N_{t+}=k+1$, corresponding with the scenario that no client leaves in the interval under study,
\[P^\uparrow_{k,u}(t):={\mathbb P}(N_{t+}=k+1, \bar B_{t+}=u+t \,|\,N_{0+}=k, \bar B_{0+}=u)=\frac{{\mathbb P}(B> u+t)}{{\mathbb P}(B> u)}.\]
Also, considering the scenario $N_{t+}=1$, in which all clients leave the system before time $t$,
\begin{equation}\nonumber P^\downarrow_{k,u}(t):=
{\mathbb P}(N_{t+}=1, \bar B_{t+}=0 \,|\,N_{0+}=k, \tau_{0+}=u)=\sum_{z=1}^{K+1}\gamma_z(u) \,{\mathbb P}\big(N_{t-}=0\,|\,N_{0+}=k, Z_{0+}=z\big),\end{equation}
where it is noticed that we have already evaluated the right-hand side of the previous display above. 
We are left with the evaluation, for $v\in(0,t)$ and $\ell = 2,\ldots,k$, of 
\begin{equation}\label{PP}p_{k\ell, uv}(t):=\frac{{\rm d}}{{\rm d}v}{\mathbb P}(N_{t+}=\ell, \bar B_{t+}\leqslant v \,|\,N_{0+}=k, \bar B_{0+}=u).\end{equation}
In order to analyze $p_{k\ell, uv}(t)$, in the sequel we work extensively with the \tredd{term}, for $v\in(0,t)$ and $\ell = 2,\ldots,k$,  \[\psi_{vt}[k,\ell]:={\mathbb P}(t-v \leqslant {\rm E}(k,\mu) \leqslant t,  {\rm E}(k,\mu)+{\rm E}(\ell-k,\mu)>t),\]
where the random variables ${\rm E}(k,\mu)$ and ${\rm E}(\ell-k,\mu)$ are independent. 
\tredd{The following lemma gives an elegant expression for these terms}.

\begin{lemma} \label{lem:3} For $v\leqslant t$ and $k<\ell$,
\[
\psi_{vt}[k,\ell] = \sum_{j=k}^{\ell-1}
\mathbb{P}({\rm Pois}(\mu t) = j)\,\mathbb{P}({\rm Bin}(j,v/t) > j-k).
\]
\end{lemma}

What is stated in Lemma \ref{lem:3} has an appealing intuitive explanation. Observe that the number of events before $t$ should be at least $k$ and should be smaller than $\ell$. This leads to the Poisson probability. Conditional on this number being $j$, more than $j-k$ of them should happen between $t-v$ and $t$, which, for each of them, happens with probability $v/t$, and in addition the corresponding events are independent. This leads to the binomial probability. A formal proof is given in \tredd{Appendix \ref{app:proofs}}.

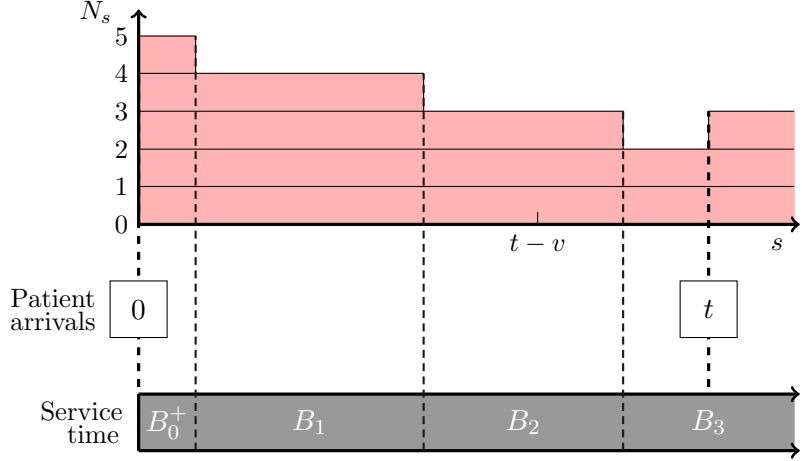
\begin{figure}[h!]
\centering
\begin{tikzpicture}[scale=0.75]

\node at (-.75, 7.75) {\small $N_s$};

\draw[fill, red!30] (0, 4) rectangle (11.5, 4.67);
\draw[fill, red!30] (0, 4.67) rectangle (11.5, 5.33);
\draw[fill, red!30] (0, 5.33) rectangle (8.5, 6);
\draw[fill, red!30] (0, 6) rectangle (5, 6.67);
\draw[fill, red!30] (0, 6.67) rectangle (1, 7.33);
\draw[fill, red!30] (10, 5.33) rectangle (11.5, 6);

\draw (0, 4.67) -- (11.5, 4.67);
\draw (0, 5.33) -- (11.5, 5.33);
\draw (0, 6) -- (8.5, 6) -- (8.5, 5.33);
\draw (0, 6.67) -- (5, 6.67) -- (5, 6);
\draw (0, 7.33) -- (1, 7.33) -- (1, 6.67);
\draw (10, 5.33) -- (10, 6) -- (11.5, 6);

\draw[very thick, ->] (0, 4) -- (0, 7.8);
\draw [very thick, ->] (0, 4) -- (11.6, 4);

\node [left] at (0, 4) {\small $0$};
\node [left] at (0, 4.67) {\small $1$};
\node [left] at (0, 5.33) {\small $2$};
\node [left] at (0, 6) {\small $3$};
\node [left] at (0, 6.67) {\small $4$};
\node [left] at (0, 7.33) {\small $5$};
\node at (11.2, 3.65) {\small $s$};

\draw[very thick, dashed] (10, 0) -- (10, 6);
\draw[very thick, dashed] (0, 0) -- (0, 6.67);
\draw[fill, black] (7, 4) -- (7, 4.22);
\node at (7, 3.65) {\small $t-v$};


\node at (-1.35, 2.7) {\small Client};
\node at (-1.5, 2.3) {\small arrivals};

\draw [fill, white] (-0.5, 2) rectangle (0.5, 3);
\draw (-0.5, 2) rectangle (0.5, 3);
\node at (0, 2.5) {$0$};
\draw[fill, white] (9.5, 2) rectangle (10.5, 3);
\draw (9.5, 2) rectangle (10.5, 3);
\node at (10, 2.5) {$t$};


\node at (-1.05, 0.7) {\small Service};
\node at (-0.8, 0.3) {\small time};

\draw[fill, gray!80] (0, 0) rectangle (11.5, 1);

\node[white] at (0.5, 0.5) {$B_0^+$};
\node[white] at (3, 0.5) {$B_1$};
\node[white] at (6.75, 0.5) {$B_2$};
\node[white] at (10, 0.5) {$B_3$};

\draw[very thick] (0, 0) -- (0, 1);
\draw[very thick, ->] (0, 1) -- (11.6, 1);
\draw[very thick, ->] (0, 0) -- (11.6, 0);

\draw[thick, densely dashed][fill, black!90] (1, 0) -- (1, 7.33);
\draw[thick, densely dashed][fill, black!90] (5, 0) -- (5, 6.67);
\draw[thick, densely dashed][fill, black!90] (8.5, 0) -- (8.5, 6);

\end{tikzpicture}
\caption{\label{fig:decomp}Example with $k=5$ and $\ell=3$. $B_0^+$ is the residual service time of client $0$ (i.e., the client in service at time $0$), and $B_i$ is the service time of client $i$.}
\end{figure}

In order to compute $p_{k\ell, uv}(t)$, we first consider the probability of $\{N_{t+}=\ell, \bar B_{t+}\leqslant v\}$ conditional on $\{N_{0+}=k, Z_{0+}=z\}$. By applying the usual procedure, we will later use the resulting probabilities to express their analogs in which the conditioning is on  $\{N_{0+}=k, \bar B_{0+}=u\}$.
For $z=1,\ldots,K$, with $I_{k\ell m z}:=(k-\ell+1)K-z+m+1$, the probability 
$q_{k\ell, z, v}(t):={\mathbb P}(N_{t+}=\ell, \bar B_{t+}\leqslant v\,|\,N_{0+}=k, Z_{0+}=z)$
equals 
 \begin{align*}\sum_{m=0}^{k-\ell+1}&{\mathbb P}({\rm Bin}(k-\ell+1,1-p)=m)
 \Big(p\,\psi_{vt}[I_{k\ell m z},I_{k\ell m z}+K] +
 (1-p)
\psi_{vt}[I_{k\ell m z},I_{k\ell m z}+K+1]\Big).\end{align*}
The reasoning behind this expression is the following.  It is first observed that we can rewrite the event $\{N_{t+} = \ell, \bar B_{t+} \leqslant v\}$ as $\{N_{(t-v){-}} \geqslant \ell, N_{t{-}} = \ell-1\}$: at time $t-$ there should be precisely $\ell-1$ clients, whereas at time $(t-v)-$ this number of clients was not reached yet. For an example we refer to Figure \ref{fig:decomp}. Hence, at time $(t-v)-$, the number of phases that have been completed is less than $K-z+1$ (certain phases of the client in service at time $0$), plus $(k-\ell)K$ (certain phases of the $k-\ell$ clients who must have left before time $t$; the last of which leaves at time $t-v$ at the latest), plus ${\rm Bin}(k-\ell+1,1-p)$ (uncertain phases of these $k-\ell+1$ clients). Adding up these numbers, one finds $I_{k\ell m z}$ in case ${\rm Bin}(k-\ell+1,1-p)=m$. At the same time, it should {\it not} be the case that by time $t$ another client has been served, i.e., the number of phases completed cannot exceed $I_{k\ell mz} + K$ if this client has $K$ phases (this occurs with probability $p$), and $I_{k\ell mz} + K+1$ if the client has $K+1$ phases. Recalling the definition of $\psi_{vt}[k,\ell]$, we thus find the probability $p\,\psi_{vt}[I_{k\ell m z},I_{k\ell m z}+K] + (1-p)\psi_{vt}[I_{k\ell m z},I_{k\ell m z}+K+1]$, as desired.

For $z=K+1$, with $I_{k\ell m}:= (k-\ell)K+m+1$, it follows analogously that $q_{k\ell, z, v}(t)$ equals
\[\sum_{m=0}^{k-\ell}{\mathbb P}({\rm Bin}(k-\ell,1-p)=m)
\Big(p\,\psi_{vt}[I_{k\ell m},I_{k\ell m}+K]
+(1-p)\psi_{vt}[I_{k\ell m},I_{k\ell m}+K+1]\Big).
\]

Combining the above expressions, we find that  \eqref{PP} equals, with $q'_{k\ell, z, v}(t)$ denoting the derivative of $q_{k\ell, z, v}(t)$ with respect to $v$,
\[p_{k\ell, uv}(t)=\sum_{z=1}^{K+1}\gamma_z(u)\,q'_{k\ell, z, v}(t).\]

 \subsubsection*{$\rhd$ Dynamic program}
 As before, we let $C_i(k,u)$,  with $i=1,\ldots,n$ and $k=1,\ldots,i$, correspond \tredd{to} the cost incurred from the arrival of the $i$-th client, given there are $k$ clients in the system immediately after the arrival of this $i$-th client and that in addition, the elapsed duration of the job in service immediately after the arrival of this $i$-th client is $u$.

 \begin{theorem} \label{thm:erl} We can determine the $C_i(k,u)$ recursively: for $i = 1,\ldots,n-1$, $k = 1,\ldots,i$, and $u\geqslant 0$,
 \begin{align*}C_i(k,u) 
 = \inf_{t\geqslant 0} \Bigg(\omega \,\bar f^\circ_{k,u}(t) &+ (1-\omega) \,\tredd{\bar h^\circ_{k,u}}
 +\sum_{\ell=2}^{k} \int_{(0,t)} p_{k\ell,uv}(t) \,C_{i+1}(\ell,v) \,{\rm d}v\\
 &+\,P^\downarrow_{k,u}(t)\, C_{i+1}(1,0) + P^\uparrow_{k,u}(t) \,C_{i+1}(k+1, u+t)\Bigg),
 \end{align*}
 whereas, for $k=1,\ldots,n$ and $u\geqslant 0$,
 \[C_n(k,u) = (1-\omega) \,\tredd{\bar h^\circ_{k,u}}.\]
 \end{theorem}


From a numerical perspective, solving the above dynamic programming problem is challenging, in particular, due to the fact that the argument $u$ can attain in principle all non-negative real values (smaller than the current time, that is). More precisely, it is a dynamic programming problem in which one of the arguments is real-valued. In practice, this means that one has to impose some sort of discretization on the variable $u$ in order to approximate the optimal strategy. The approach we propose here is to evaluate, for a given $\Delta>0$, $\xi_i(k,m)$, which is to be interpreted as an  approximation for $C_i(k,m\Delta)$ that becomes increasingly accurate as $\Delta\downarrow 0$, for  $k = 1,\ldots,i$ and $m\in{\mathbb N}_0$ \tred{(with ${\mathbb N}_0:={\mathbb N}\cup\{0\}=\{0,1,\ldots\}$).}
Define by $[x]_y$ the truncation of $x$ to the interval $[0,y]$, for some $y>0$:
\[[x]_y := \min\{\max\{0,x\}, y\}.\]
Also, for  $t\in {\mathbb N}_0$ and $\Delta>0$,
\begin{align*}
\bar q_{k\ell, mj}(t)&:=\int_{[(j-\frac{1}{2})\Delta]_{t\Delta}}^{[(j+\frac{1}{2})\Delta]_{t\Delta}} p_{k\ell, m\Delta,v}(t\Delta){\rm d}v\\&=
\sum_{z=1}^{K+1} \gamma_z(m\Delta) \Bigg(q_{k\ell,z,[(j+\frac{1}{2})\Delta]_{t\Delta}}(t\Delta)-
q_{k\ell,z,[(j-\frac{1}{2})\Delta]_{t\Delta}}(t\Delta)
\Bigg).\end{align*}
This quantity is to be interpreted as an approximation of the transition probability of the number of clients jointly with the age of the client in service (truncated to a multiple of $\Delta$, that is), from the state $(k,m\Delta)$ to the state $(\ell,j\Delta)$, over a time interval of length $t\Delta$, for some $t\in {\mathbb N}_0$.
The approximation of the dynamic programming recursion thus becomes
\begin{align}\xi_i(k,m)
 \label{eqn:xi_ikm}
 \nonumber
 = \inf_{t\in{\mathbb N}_0} \Bigg(\omega \,\bar f^\circ_{k,m\Delta}(t\Delta) &+ (1-\omega) \,\tredd{\bar h^\circ_{k,m\Delta}}
 +\sum_{\ell=2}^{k} \sum_{j=0}^{t} \bar q_{k\ell, mj}(t)\,\xi_{i+1}(\ell,j)\\
 +\,P^\downarrow_{k,m\Delta}(t\Delta)\, &\xi_{i+1}(1,0) + P^\uparrow_{k,m\Delta}(t\Delta) \,\xi_{i+1}(k+1, m+t)\Bigg),
 \end{align}
 whereas, for $k=1,\ldots,n$ and $m\in{\mathbb N}_0$,
 \[\xi_n(k,m) = (1-\omega) \tredd{\,\bar h^\circ_{k,m\Delta}}.\]
 
\subsection{Hyperexponential distribution} \label{HED} In this case the service time $B$ equals with probability $p\in[0,1]$ an exponentially distributed random variable with mean $\mu_1^{-1}$, and with probability $1-p$ an exponentially distributed random variable with mean $\mu_2^{-1}$. This means that 
\[T = \left(\begin{array}{cc}
-\mu_1&0\\
0&-\mu_2\end{array}\right),\]
and $\alpha_1 = p$\tredd{, i.e., ${\bs \alpha} = (p,1-p)^T$.} 

As in the weighted Erlang case, the phase is not observable: it cannot be observed which of the two exponential random variables $B$ corresponds to. As before, we determine the distribution of the phase conditional on the elapsed service time. It is checked easily
that $\gamma_z(t) =  \gamma^\circ_z(t)/\gamma^\circ(t)$, where
\begin{align*} \gamma^\circ(t)&:={\mathbb P}(B>t)
=p \,e^{-\mu_1 t}+(1-p)\,e^{-\mu_2 t},\\
 \gamma^\circ_z(t)&:={\mathbb P}(Z_t=z, B>t)= p \,e^{-\mu_1 t}\mathbbm{1}_{\{z=1\}}+(1-p)\,e^{-\mu_2 t}\mathbbm{1}_{\{z=2\}}.
 \end{align*}
 As a sanity check, observe that $\gamma_1(0)=1-\gamma_2(0)=p$ and (recalling that we imposed $\mu_1>\mu_2$) $\gamma_1(\infty)=1-\gamma_2(\infty)=0$.

Mimicking the structure we have used in the weighted Erlang case, we subsequently discuss the mean idle time, the mean waiting time, the transition probabilities, and the resulting dynamic program. 

 \subsubsection*{$\rhd$ Mean idle time}
Again we first determine the mean idle time, given $Z_{0+}= z$:
 \begin{align*}\bar f_{k z }(t)&:= \int_0^t {\mathbb P}\big(N_s=0\,|\, N_{0+}=k ,Z_{0+}= z \big)\,{\rm d}s.\end{align*}
 Translating the expression for $\bar f_{k z }(t)$ into an expression for $\bar f^\circ_{k,u }(t)$ (i.e., a probability in which the conditioning is on information that is observable) works as in the weighted Erlang case:
\begin{equation}\label{f_ku_circ_hyp}
     \tredd{\bar{f}^{\circ}_{k,u}(t) := \int_{0}^{t}\mathbb{P}(N_s = 0\,|\, N_{0+} = k, \bar{B}_{0+} = u)\,{\rm d}s = \gamma_1(u)\bar{f}_{k1}(t) + \gamma_2(u)\bar{f}_{k2}(t).}
\end{equation}

 Below we concentrate on the case $z=1$, but evidently the case $z=2$ can be dealt with fully analogously due to the inherent symmetry of the hyperexponential distribution.  
 The starting point is the identity 
 \[\bar f_{k 1}(t) = 
 \sum_{m=0}^{k-1}\mathbb{P}(\text{Bin}(k-1,p) = m) \sigma_t[m+1, k-1-m],\]
 where, with ${\rm E}(m,\mu_1)$ and ${\rm E}(k,\mu_2)$ being independent,
 \[\sigma_t[m,k]:= \int_0^t{\mathbb P}({\rm E}(m,\mu_1) + {\rm E}(k,\mu_2) \leqslant s)\,{\rm d}s.\]
 The following lemma is proved in \tredd{Appendix \ref{app:proofs}}.
 \begin{lemma}\label{L4}
 For $m,k\in\{1,2,\ldots\}$ and $t\geqslant 0$,
\[ \sigma_t[m,k] = \left(t-\frac{k}{\mu_2}\right){\mathbb P}\left(
   {\rm E}(m,\mu_1)\leqslant t\right)
             - \frac{m}{\mu_1}{\mathbb P}\left(
   {\rm E}(m+1,\mu_1)\leqslant t\right)+\frac{\mu_1}{\mu_2}\sum_{i=0}^{k-1}(k-i) 
       \rho_t[m-1,i]       ,
     \]
   where
     \[\rho_t[m,k]:= \int_0^t r_{t,u}[m,k] \,{\rm d}u,\:\:\:\:
     r_{t,u}[m,k] 
     :=
     e^{-\mu_1 u}\frac{(\mu_1 u)^m}{m!} e^{-\mu_2 (t-u)}\frac{(\mu_2 (t-u))^k}{k!}{\rm d}u.\] 
 \end{lemma}

   We proceed by pointing out how the \tredd{term} $\rho_t[m,k]$ can be evaluated. To this end, we set up a  recursive procedure.   Integration by parts yields that, for $k,m=1,2,\ldots,$
  \[\rho_{t}[m,k] = \rho_{t}[m-1,k]+\frac{\mu_2}{\mu_1}\rho_{t}[m,k]-\frac{\mu_2}{\mu_1}\rho_{t}[m,k-1],\]
  which after rearranging yields the recursive relation
   \[\rho_{t}[m,k]  = \frac{\mu_1  \rho_t[m-1,k] -\mu_2\rho_t[m,k-1]}{\mu_1-\mu_2}.\]
In addition, $\rho_{t}[m,0]$ and $\rho_{t}[0,k]$ can be given in closed form. Indeed (with  $\mu_1\not=\mu_2$),
\begin{align*}\rho_{t}[m,0] &= \frac{\mu_1^m\,e^{-\mu_2 t}}{(\mu_1-\mu_2)^{m+1}}\,\left(1-\sum_{i=0}^m e^{(\mu_2-\mu_1)t}\frac{((\mu_1-\mu_2)t)^{i}}{i!}\right),\\
\rho_{t}[0,k] &= \frac{\mu_2^k\,e^{-\mu_1 t}}{(\mu_2-\mu_1)^{k+1}}\,\left(
1
-
\sum_{i=0}^k e^{(\mu_1-\mu_2)t}\frac{((\mu_2-\mu_1)t)^{i}}{i!}
\right).\end{align*}
Using the above recursive relation in combination with these initial values, we can evaluate the quantities  $\rho_{t}[m,k]$.

\subsubsection*{$\rhd$ Mean waiting time}
\tredd{Again the evaluation of $\bar g^{\circ}_{k,u}(t)$ is provided in Appendix \ref{app:mean_idle_time}, to proceed with $\bar h^{\circ}_{k,u}$ directly.} 
It is not hard to see that $\bar h_{1z}=0$, whereas for $k=2,\ldots,i$,
 \[\bar h_{kz} :=\frac{1}{\mu_1}\,\mathbbm{1}_{\{z=1\}}+ \frac{1}{\mu_2}\mathbbm{1}_{\{z=2\}} +\frac{k-2}{\mu},\]
where the mean waiting time equals ${1}/{\mu} = {p}/{\mu_1} + {(1-p)}/{\mu_2}$. \tredd{We thus obtain}
\begin{equation}\label{h_ku_circ_hyp}
    \tredd{\bar{h}_{k,u}^{\circ} := \gamma_1(u)\bar{h}_{k1} + \gamma_2(u)\bar{h}_{k2}.}
\end{equation}

 \subsubsection*{$\rhd$ Transition probabilities}
To this end, we continue by analyzing the probability $q_{k\ell, z, v}(t):={\mathbb P}(N_{t+}=\ell, \bar B_{t+}\leqslant v\,|\,N_{0+}=k, Z_{0+}=z)$. We start by considering the case that $\ell\in\{2,\ldots,k\}$, i.e., for the moment we exclude the scenarios in which no or all clients leave before time $t$. As we have observed, we can rewrite $\{N_{t+} = \ell, \bar B_{t+} \leqslant v\}$ as $\{N_{(t-v){-}} \geqslant \ell, N_{t{-}} = \ell-1\}$. Define, for $i = 1,2$,
 \[\chi_{vt,i}[k,\ell]:= {\mathbb P}(t-v \leqslant {\rm E}(k,\mu_1)+{\rm E}(\ell,\mu_2)\leqslant t,
 {\rm E}(k,\mu_1)+{\rm E}(\ell,\mu_2) + {\rm E}(1,\mu_i)>t),\]
 with $ {\rm E}(k,\mu_1)$, ${\rm E}(\ell,\mu_2)$, and ${\rm E}(1,\mu_i)$ independent \tredd{Erlang distributions}. The following lemma, proved in \tredd{Appendix \ref{app:proofs}}, provides expressions for $\chi_{vt,1}[k,\ell]$ and $\chi_{vt,2}[k,\ell]$.
 
 \begin{lemma} \label{L5} For $v\leqslant t$ and $k,\ell=1,2,\ldots$,
 \begin{align*}\chi_{vt,1}[k,\ell] &= \mu_2\,\rho_{t}[k,\ell-1]-e^{-\mu_1 v}\mu_2 \,\rho_{t-v}[k,\ell-1],\\
\chi_{vt,2}[k,\ell] &=\mu_1\,\rho_{t}[k-1,\ell]-e^{-\mu_2 v}\mu_1 \,\rho_{t-v}[k-1,\ell].
  \end{align*}
 \end{lemma}

It is relatively straightforward to provide closed-form expressions for the probabilities $\chi_{vt,1}[k,\ell]$ and $\chi_{vt,2}[k,\ell]$ in case $k$ and/or $\ell$ equals $0$, thus complementing the cases dealt with in Lemma \ref{L5}. Indeed, for $k,\ell=1,2,\ldots$,
\begin{align*}
\chi_{vt,1}&[0,\ell]=
\int_{t-v}^t \int_{t-s}^\infty \mu_2 e^{-\mu_2 s} \frac{(\mu_2 s)^{\ell-1}}{(\ell-1)!}\,\mu_1e^{-\mu_1 u}\,{\rm d}u\,{\rm d}s\\
&=e^{-\mu_1 t}
\left(\frac{\mu_2}{\mu_2-\mu_1}\right)^\ell 
\left(\sum_{i=0}^{\ell-1}e^{-(\mu_2-\mu_1)(t-v)}\frac{((\mu_2-\mu_1)(t-v))^i}{i!}-\sum_{i=0}^{\ell-1}e^{-(\mu_2-\mu_1)t}\frac{((\mu_2-\mu_1)t)^i}{i!}\right), \\
\chi_{vt,1}&[k,0] = {\mathbb P}({\rm Pois}(\mu_1 t)=k)\,{\mathbb P}({\rm Bin}(k,v/t)>0)= e^{-\mu_1 t}\frac{\mu_1^k}{k!}\big( t^k -(t-v)^k\big),
\end{align*}
where the latter statement follows due to Lemma \ref{lem:3}. The expressions for $\chi_{vt,2}[0,\ell]$ and $\chi_{vt,2}[k,0]$ follow by symmetry: they equal $\chi_{vt,1}[\ell,0]$ and $\chi_{vt,1}[0,k]$, respectively, but with the roles of $\mu_1$ and $\mu_2$ being swapped. Also, $\chi_{vt,i}[0,0] = e^{-\mu_i t} \,\mathbbm{1}_{\{v=t\}}.$
 
 \vb
 
 We thus obtain that
 \begin{align*}q_{k\ell, 1, v}(t)&=\sum_{m=0}^{k-\ell}{\mathbb P}({\rm Bin}(k-\ell,p)=m)
\Big(p\,\chi_{vt,1}[m+1, k-\ell-m]
+(1-p)\chi_{vt,2}[m+1,k-\ell-m]\Big),\\
q_{k\ell, 2, v}(t)&=\sum_{m=0}^{k-\ell}{\mathbb P}({\rm Bin}(k-\ell,p)=m)
\Big(p\,\chi_{vt,1}[m, k-\ell-m+1]
+(1-p)\chi_{vt,2}[m,k-\ell-m+1]\Big).
\end{align*}
\tredd{Now, with $q'_{k\ell,z,v}(t)$ denoting the derivative of $q_{k\ell,z,v}(t)$ with respect to $v$,}
\begin{equation}\label{p_kl_uv_hyp}
    \tredd{p_{k\ell,u,v}(t) = \gamma_1(u)\,q'_{k\ell,1,v}(t) + \gamma_2(u)\,q'_{k\ell,2,v}(t).}
\end{equation}
We are left with analyzing the cases $N_{t+}=k+1$ and $N_{t+}=1$.
Recall that $N_{t+}=k+1$ corresponds with the scenario that no client leaves in the interval under study,
\begin{equation}\label{P_up_hyp}
\tredd{P^\uparrow_{k,u}(t):={\mathbb P}(N_{t+}=k+1, \bar B_{t+}=u+t \,|\,N_{0+}=k, \bar B_{0+}=u)=\frac{{\mathbb P}(B> u+t)}{{\mathbb P}(B> u)}.}
\end{equation}
Regarding the scenario $N_{t+}=1$, with  all clients leaving  before time $t$,
{\small
\begin{equation}\label{P_down_hyp}
\tredd{P^\downarrow_{k,u}(t):=
{\mathbb P}(N_{t+}=1, \bar B_{t+}=0 \,|\,N_{0+}=k, \bar{B}_{0+}=u)=\sum_{z=1}^{2}\gamma_z(u) \,{\mathbb P}\big(N_{t-}=0\,|\,N_{0+}=k, Z_{0+}=z\big).}
\end{equation}}\noindent
Here it is noticed that the probabilities appearing on the right-hand side can be evaluated with arguments similar to the ones used above (in particular, use that the number of clients waiting at time $0$ corresponding with an exponential distribution with rate $\mu_1$ being binomially distributed with parameters $k-1$ and $p$). 

\subsubsection*{$\rhd$ Dynamic program}
We are now in a position to define our dynamic programming algorithm. The cost $C_i(k,u)$ is as defined before. 
 
 \begin{theorem} We can determine the $C_i(k,u)$ recursively with the procedure of Theorem \ref{thm:erl}, but with the \tredd{$\bar f^\circ_{k,u}(t)$, $\bar h_{k,u}^\circ$, $p_{k\ell,u,v}(t)$, $P^\uparrow_{k,u}(t)$, and $P^\downarrow_{k,u}(t)$ replaced by their counterparts for the hyperexponential case, i.e., equations \eqref{f_ku_circ_hyp}, \eqref{h_ku_circ_hyp}, \eqref{p_kl_uv_hyp}, \eqref{P_up_hyp} and \eqref{P_down_hyp}, respectively.}
  \end{theorem}
  
  As in the case of the weighted Erlang service times, this procedure can be numerically evaluated by discretizing the elapsed service time $u$ (as multiples of $\Delta>0$, that is), and by searching over arrival times that are also multiples of this $\Delta$.
  
  \section{Numerical evaluation}\label{num}
  As mentioned above, primarily due to the fact that we have to compute the optimal $t\geqslant 0$ (or $t\in{\mathbb N}_0$, in the discretized version) for any value of the elapsed service time $u$, performing the dynamic programming procedure is numerically challenging. 
  In this section we present our applet that provides dynamic appointment schedules, provide more detail about its implementation, and discuss a number of illustrative experiments.
  
  \subsection{Applet} We have developed an applet providing the optimal appointment time of the next client for a given instance.\tredd{\footnote{The applet can be approached through
  {\tt \url{https://dynamicschedule.eu.pythonanywhere.com}.}}} In this context, an {\it instance} is a combination of the SCV of the \tred{(homogeneous)} service times ${\mathbb S}(B)$, the weight $\omega$, the number of clients $n$, the index $i\in\{1,\ldots,n\}$ of the entering client, the number of clients in the system $k\in\{1,\ldots,i\}$ when the $i$-th client enters (including this entering client), and the elapsed service time $u\geqslant 0$ of the client in service (if any). To cover situations we came across in the service system literature, we let ${\mathbb S}(B)\in [0.2, 2]$. The weight $\omega$ can be chosen in the set $\omega \in\{0.1,0.2,\ldots, 0.9\}$, and one can choose the number of clients $n \in \{1,2,\dots,20\}$. 
  Note that above we normalized time such that ${\mathbb E}B=1$, whereas in the applet ${\mathbb E}B$ can be freely chosen.
  
  
 \subsection{Implementation}

To avoid doing computations in real-time, we have chosen an approach in which we have off-line performed the computations on a dense grid of instances, varying $\omega$ and the SCV. 
\tred{It is inevitable to follow an approach that uses pre-computed values, due to the intrinsic complexity of  our dynamic programming problem, involving multiple dimensions of which one (the elapsed service time) is effectively continuous; this can be considered as a manifestation of the notorious curse of dimensionality.}

It is important to note that the majority of the computation time of $\xi_i(k,m)$, for a given instance, lies in finding the value of $t$ that minimizes the expression in \eqref{eqn:xi_ikm}, which we denote by $\tau_i(k, m)$. Starting with $i=n$ and decreasing $i$ in steps of 1, we have gained a significant speed-up by storing all computed values of $\xi_{i+1}(\cdot, \cdot)$ and caching the $\bar q_{k\ell, mj}(t)$ to avoid multiple computations of the same quantity. To quickly find the optimal value of $t$, it is essential to have a good initial guess. We have employed the following three-step approach for each value of $\omega$:
\begin{enumerate}
\item Observe that ${\mathbb S}(B)=1$ corresponds to the relatively straightforward case of exponentially distributed service times, as dealt with in Section~\ref{homexp}. Hence, Step~1 consists of computing $\xi_i(k, m)$ and $\tau_i(k,m)$ for ${\mathbb S}(B)=1$, for $k=1,\dots,i$. Note that $\xi_i(k, m)$ and $\tau_i(k,m)$ do not depend on $m$ due to the memoryless property of the exponential distribution.
\item In Step 2, our goal is to obtain good initial estimates of $\tau_i(k,m)$. We choose a relatively large step size $\Delta$, which reduces the computation time of \eqref{eqn:xi_ikm}, to obtain good initial estimates for $\tau_i(k,m)$. We compute for ${\mathbb S}(B)\neq 1$ the values of $\tau_i(k,m)$ only for $i=n-1$ and $i=n-2$ and $m=0$. The reason \tredd{for computing} only these values, is that $\tau_{n-2}(k, 0)$ is an excellent approximation for $\tau_{i}(k, 0)$ for $i=1,\dots,n-3$. The value of $\tau_{n-1}(k, m)$ is different because in this case, only the last client is suffering from potentially high waiting times.
\item In the final step, we take a smaller value of $\Delta$ (typically a factor 10 smaller) to compute $\xi_i(k, m)$ and $\tau_i(k, m)$ for all values of $i$, $k$, and $m$. Finding the minimum of \eqref{eqn:xi_ikm} is done by a nearest neighbor search starting in the estimate obtained in Step~2 and stopping when the absolute decrease in the objective function is smaller than $\epsilon=10^{-5}$. Since $\Delta$ is small, the difference between $\xi_i(k, m)$ and $\xi_i(k, m+1)$ is very small as well. For this reason, we have decided to compute $\xi_i(k, m)$ only for $m\in \{0,10,20,\dots, 250\}$ and obtain the other values by interpolation or extrapolation.
\end{enumerate}     

A final remark: when computing the optimal schedule (i.e., the arrival times) for a given value of $n$, say $n=N$, we can automatically derive the schedules for all $n=1,2,\dots, N-1$ from this schedule. To clarify this, suppose that we have $n=N$ \tredd{clients} that need to be scheduled. Our algorithm computes the optimal arrival times $\tau_i^{(N)}(k, m)$ and the corresponding costs $\xi_i^{(N)}(k, m)$ for $i=1,\dots,N$, where we use the superscript $(N)$ to indicate that these quantities correspond to the schedule for $N$ \tredd{clients}. Using this schedule, we can immediately obtain the optimal schedule for scheduling $N' < N$ \tredd{clients}. Indeed, due to the assumption of homogeneous service times, the transition probabilities and the expected idle and waiting times do not depend on the client position $i$, and therefore
\[
\xi_i^{(N')}(k, m) = \xi_{i+N-N'}^{(N)}(k, m) \quad\text{ and }\quad
\tau_i^{(N')}(k, m) = \tau_{i+N-N'}^{(N)}(k, m).
\]



  \subsection{Experiments}
  \label{subsec:experiments}
  In this subsection we present a series of experiments that highlight the performance of our dynamic scheduling approach. We start by providing insight into the efficiency gain, to later study in greater detail the role played by the model parameters.
  
  \subsubsection*{$\rhd$~Experiment 1: Efficiency gain} 
\begin{figure}
      \centering
      \includegraphics[width=0.48\textwidth]{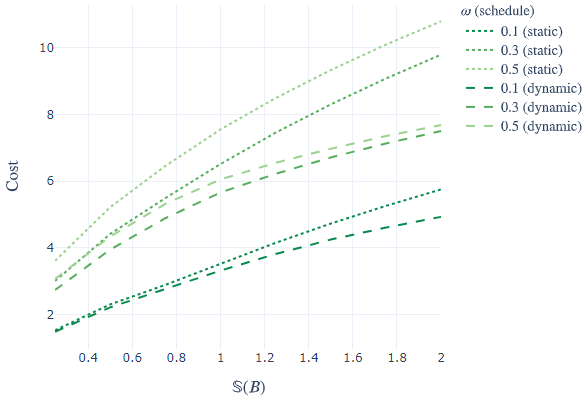}\:
      \includegraphics[width=0.48\textwidth]{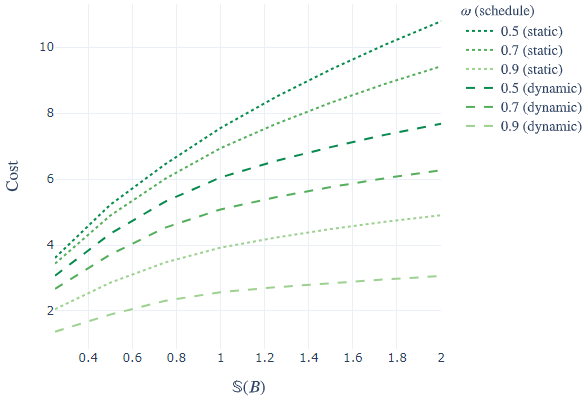}
            \caption{Cost of precalculated and dynamic schedule as function of ${\mathbb S}(B)$ for different values of $\omega$. Left panel: $\omega\in\{0.1,0.3,0.5\}$. Right panel: $\omega\in\{0.5,0.7,0.9\}$.}
            \label{fig:gain}
  \end{figure}

The first experiment quantifies the gain of using dynamic schedules relative to using precalculated schedules. In Figure \ref{fig:gain} we show for $\omega\in\{0.1,0.3,\dots,0.9\}$ and $n = 15$, as a function of ${\mathbb S}(B)$, the cost of the dynamic schedule and the cost of the precalculated schedule. In Table \ref{adapt_vs_prec_phase}, we see the ratio between the two costs. We conclude that using dynamic schedules particularly pays off for high values of the weight (i.e., $\omega$), \tred{in line with what we discussed in Remark \ref{RRR}}, and high values of the SCV of the service times (i.e., ${\mathbb S}(B)$).

  {\small
  \begin{table}
  \centering
  \begin{tabular}{|c|c|ccccccccc|}
  \hline
  $\mathbb{S}(B)$ & $\omega$                              & 0.1  & 0.2  & 0.3  & 0.4  & 0.5   & 0.6  & 0.7  & 0.8  & 0.9  \\
  \hline
  0.25            & $K_{\text{dyn}}(\mathbb{S}(B),\omega)$ & 1.49 & 2.26 & 2.74 & 3.00 & 3.07  & 2.97 & 2.67 & 2.16 & 1.37 \\
                  & $K_{\text{pre}}(\mathbb{S}(B),\omega)$ & 1.53 & 2.41 & 3.01 & 3.40 & 3.61  & 3.63 & 3.44 & 2.96 & 2.06 \\
                  & $r(\mathbb{S}(B),\omega)$             & 0.97 & 0.94 & 0.91 & 0.88 & 0.85  & 0.82 & 0.78 & 0.73 & 0.67 \\
  \hline
  0.5             & $K_{\text{dyn}}(\mathbb{S}(B),\omega)$ & 2.22 & 3.31 & 3.95 & 4.28 & 4.34  & 4.15 & 3.71 & 2.99 & 1.89 \\
                  & $K_{\text{pre}}(\mathbb{S}(B),\omega)$ & 2.31 & 3.57 & 4.42 & 4.96 & 5.22  & 5.21 & 4.89 & 4.18 & 2.86 \\
                  & $r(\mathbb{S}(B),\omega)$             & 0.96 & 0.93 & 0.89 & 0.86 & 0.83  & 0.80 & 0.76 & 0.72 & 0.66 \\
  \hline
  0.75            & $K_{\text{dyn}}(\mathbb{S}(B),\omega)$ & 2.77 & 4.11 & 4.89 & 5.27 & 5.32  & 5.07 & 4.53 & 3.64 & 2.31 \\
                  & $K_{\text{pre}}(\mathbb{S}(B),\omega)$ & 2.89 & 4.46 & 5.49 & 6.14 & 6.45  & 6.42 & 6.01 & 5.11 & 3.47 \\
                  & $r(\mathbb{S}(B),\omega)$             & 0.96 & 0.92 & 0.89 & 0.86 & 0.82  & 0.79 & 0.75 & 0.71 & 0.66 \\
  \hline
  1               & $K_{\text{dyn}}(\mathbb{S}(B),\omega)$ & 3.32 & 4.83 & 5.66 & 6.04 & 6.05  & 5.73 & 5.08 & 4.07 & 2.57 \\
                  & $K_{\text{pre}}(\mathbb{S}(B),\omega)$ & 3.51 & 5.33 & 6.51 & 7.23 & 7.55  & 7.47 & 6.94 & 5.85 & 3.92 \\
                  & $r(\mathbb{S}(B),\omega)$             & 0.95 & 0.91 & 0.87 & 0.83 & 0.80  & 0.77 & 0.73 & 0.70 & 0.66 \\
  \hline
  1.25            & $K_{\text{dyn}}(\mathbb{S}(B),\omega)$ & 3.82 & 5.39 & 6.22 & 6.57 & 6.55  & 6.17 & 5.45 & 4.33 & 2.72 \\
                  & $K_{\text{pre}}(\mathbb{S}(B),\omega)$ & 4.15 & 6.18 & 7.45 & 8.20 & 8.49  & 8.33 & 7.67 & 6.40 & 4.23 \\
                  & $r(\mathbb{S}(B),\omega)$             & 0.92 & 0.87 & 0.83 & 0.80 & 0.77  & 0.74 & 0.71 & 0.68 & 0.64 \\
  \hline
  1.5             & $K_{\text{dyn}}(\mathbb{S}(B),\omega)$ & 4.25 & 5.87 & 6.71 & 7.04 & 6.97  & 6.55 & 5.76 & 4.56 & 2.85 \\
                  & $K_{\text{pre}}(\mathbb{S}(B),\omega)$ & 4.73 & 6.94 & 8.30 & 9.07 & 9.33  & 9.09 & 8.32 & 6.88 & 4.49 \\
                  & $r(\mathbb{S}(B),\omega)$             & 0.90 & 0.85 & 0.81 & 0.78 & 0.75  & 0.72 & 0.69 & 0.66 & 0.64 \\
  \hline
  1.75            & $K_{\text{dyn}}(\mathbb{S}(B),\omega)$ & 4.61 & 6.29 & 7.13 & 7.44 & 7.35  & 6.88 & 6.03 & 4.76 & 2.96 \\
                  & $K_{\text{pre}}(\mathbb{S}(B),\omega)$ & 5.26 & 7.64 & 9.07 & 9.86 & 10.09 & 9.78 & 8.90 & 7.31 & 4.71 \\
                  & $r(\mathbb{S}(B),\omega)$             & 0.88 & 0.82 & 0.79 & 0.75 & 0.73  & 0.70 & 0.68 & 0.65 & 0.63 \\
  \hline
  \end{tabular}
  \caption{Cost of dynamic and precalculated schedule in Experiment~1.}
  \label{adapt_vs_prec_phase}
  \end{table}}

\tredd{
Obviously, in practice, one does not always have good a priori  estimates of (the parameters of) the service-time distributions. We continue by doing a simulation experiment where we sample a \emph{random} SCV value in each run. This will provide more insight in the impact of fluctuations in the SCVs on the quality of the dynamic schedules. Additionally, it enables us to create confidence intervals (CIs) for the mean decrease in costs. The SCV of the clients is sampled uniformly at random from the interval $[0.5, 1.5]$, while the mean service time is set equal to $1$. For each set of parameter values, we determine the optimal precalculated schedule and the optimal dynamic schedule and compare their costs. We repeat this simulation experiment 100,000 times and compute 95\% CIs for the mean difference of the costs of precalculated and dynamic scheduling, for $\omega\in\{0.2, 0.5, 0.8\}$. 
The results, shown in Table~\ref{tbl:randomparametersgeneral}, confirm that there is indeed  a significant difference. For practical applications, the prediction intervals~(PIs) are arguably even more interesting. Each prediction interval consists of the 2.5\% and the 97.5\% percentiles of the simulated costs, which is an excellent quantification of the magnitude of the random fluctuations one can expect in the costs. The PIs indicate that fluctuations of more than 30\% can be expected due to the randomness in the parameters, confirming that it is really important for practical implementation of the scheduling policies to have accurate estimates of the service-time distributions. In this example, we have varied the variance of the service-time distributions; in Example~\ref{E3} in Appendix~\ref{app:hetexp}, we study the impact of varying their means in a heterogeneous case.
}

\begin{table}[ht]
\centering
\tredd{
\begin{tabular}{|c|ccc|cc|}
\hline
 $\omega$  &  $K_{\text{dyn}}(\omega)$ & $K_{\text{pre}}(\omega)$ &
 $K_{\text{pre}}(\omega) - K_{\text{dyn}}(\omega)$ &
 95\% CI & 95\% PI\\
\hline      
0.2   & 4.738 & 5.312 & 0.574 & [0.573, 0.576] & [0.268, 1.050] \\
0.5   & 5.905 & 7.452 & 1.546 & [1.543, 1.549] & [0.894, 2.323] \\ 
0.8   & 3.963 & 5.732 & 1.769 & [1.767, 1.771] & [1.211, 2.305] \\
\hline
\end{tabular}}
\caption{\tredd{Cost of dynamic and precalculated schedule, with $n=15$ and random $\mathbb{S}(B)\in[0.5,1.5]$.  95\% confidence intervals and prediction intervals are computed for the difference in costs.}}
\label{tbl:randomparametersgeneral}
\end{table}

  \subsubsection*{$\rhd$~Experiment 2: Impact of model parameters} We continue by studying the impact of incorporating the elapsed service time $u$. As argued before, in the case of exponentially distributed service times this elapsed service time has no impact due to the memoryless property. Define by $\tau_i(k,u)$ the optimizing argument in the definition of $C_i(k,u)$.
  
  \begin{figure}
      \centering
      \includegraphics[width=0.48\textwidth]{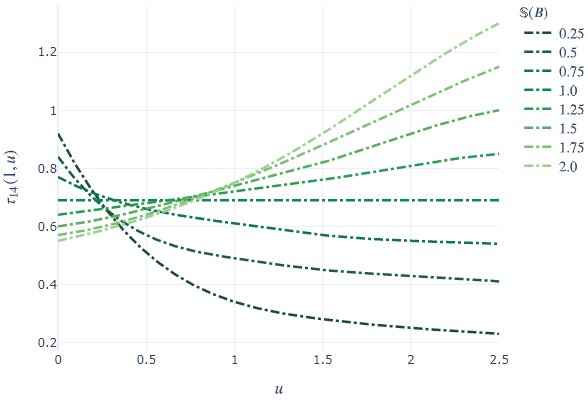}\:\:
      \includegraphics[width=0.48\textwidth]{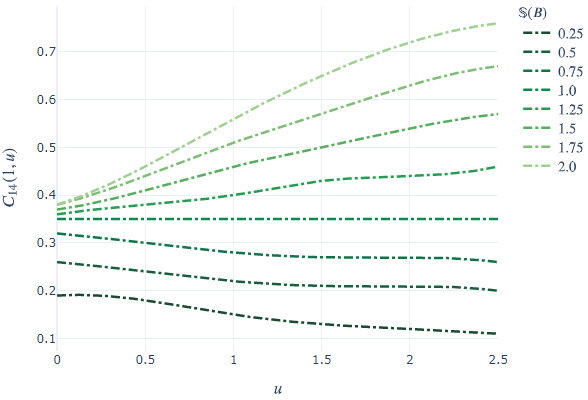}
            \caption{Optimized interarrival time $\tau_{14}(1,u)$ and corresponding cost $C_{14}(1,u)$ as a function of $u\geqslant 0$ for several values of ${\mathbb S}(B)$ and $\omega=0.5$. \label{Optu}}
  \end{figure}
  
  

  In Figure \ref{Optu} we show the dependence of the optimal interarrival time $\tau_i(k,u)$ and $C_i(k,u)$ on $u$. For illustrational purposes, we do this for specific values of $i$, $k$, and $\omega$, but the conclusions below hold in general.  In the first place, we observe that in case ${\mathbb S}(B)=1$ the elapsed service time has no impact on $\tau_i(k,u)$ and $C_i(k,u)$, in line with the reasoning above. In case ${\mathbb S}(B)<1$, the service times are `relatively deterministic'. 
  This explains why in this regime both $\tau_i(k,u)$ and $C_i(k,u)$ decrease in $u$: the longer the elapsed part of the current service, the sooner it is expected to finish, and the more accurately this completion epoch can be predicted.
  If, on the other hand, ${\mathbb S}(B)>1$, then the service times are `less deterministic' than when stemming from the exponential distribution, so that $\tau_i(k,u)$ and $C_i(k,u)$ increase in $u$. \tred{In the latter case, in which we work with the hyperexponential distribution, the remaining service time is increasing in $u$: with increasing probability the exponential distribution corresponding to the minimum of $\mu_1$ and $\mu_2$ (i.e., the one with the highest mean) is sampled from.}
  

  \subsubsection*{$\rhd$~Experiment 3: Impact of ignoring coefficient of variation}
\begin{figure}
    \centering
    \begin{floatrow}
      \ffigbox[\FBwidth]{\caption{\tredd{Ratio of cost when using a dynamic schedule based on the correct SCV and cost when using a precalculated schedule based on an SCV equal to one for $n = 15$ and several values of $\omega$.}}\label{SCV_wrong-1}}{%
        \includegraphics[width=0.46\textwidth]{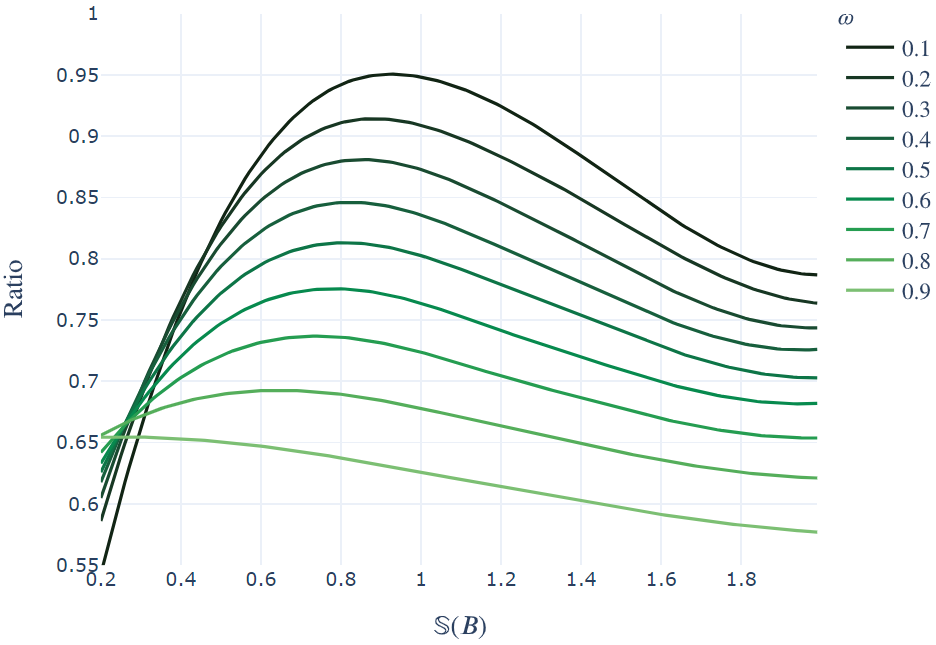}
      }
      \ffigbox[\FBwidth]{\caption{\tredd{Ratio of cost when using a dynamic schedule based on the correct SCV and cost when using a dynamic schedule based on an SCV equal to one for $n = 15$ and several values of $\omega$.}}\label{SCV_wrong-2}}{%
         \includegraphics[width=0.46\textwidth]{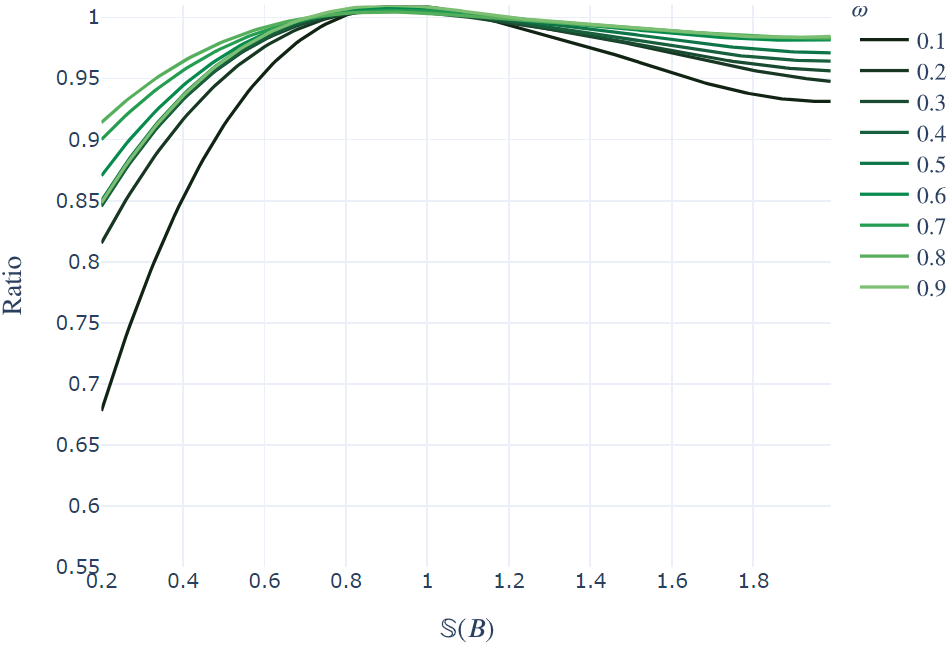}
      }
    \end{floatrow}
  \end{figure}

  In many approaches, owing to its convenient properties, it is assumed that service times have an exponential distribution. In this experiment, our objective is to quantify the loss of efficiency due to ignoring the non-exponentiality of the service-time distribution. In the experiment performed we let the service times be endowed with an SCV different from one. The following three different schedules will be compared:
  \begin{itemize}
      \item[(i)]a \emph{static} schedule assuming \emph{exponential} service times, computed using the methodology described in e.g.\ \cite{KUIP};
      \item[(ii)]a \emph{dynamic} schedule assuming \emph{exponential} service times, computed using the DP-based machinery developed in Section \ref{homexp};
      \item[(iii)]a \emph{dynamic} schedule assuming \emph{the correct} SCV, computed using our methods developed in Section~\ref{phase}.
      \end{itemize}
  As before, all mean service times are set equal to one. We estimate the cost of~(i) and~(ii) by simulation and compare it with the cost of (iii).
  As can be seen in \tredd{Figure \ref{SCV_wrong-1} and Figure \ref{SCV_wrong-2}}, wrongly assuming exponentiality may lead to a significant loss of efficiency. Below we provide more detailed observations. 
  \begin{itemize}
      \item[$\circ$]
  Particularly when comparing (i) and (iii), which we do in \tredd{Figure \ref{SCV_wrong-1}}, we observe a substantial loss: using a precalculated schedule based on SCV equal to 1 performs considerably worse than the dynamic schedule based on the correct SCV. This conclusion is in line with what was observed in Experiment~1.
  \item[$\circ$]
  When comparing (ii) and (iii), which we do in \tredd{Figure \ref{SCV_wrong-2}}, the loss is relatively small for values of the SCV around 1. This may have to do with the fact that in dynamic schedules there are ample opportunities to correct when unforeseen scenarios occur, thus mitigating the effect of the randomness in the service times. We observe, however, that for lower values of the SCV, the loss can be substantial; it is noted that in the healthcare context such low values  are common, as was reported in \cite{CV}. 
  \end{itemize}

  \section{Discussion}\label{sec:discu}
  \tred{In this section we reflect on two relevant aspects of our approach. The first is related to its {\it robustness}: what is the impact of replacing general distributions by their phase-type counterparts (\tredd{i.e., phase-type distributions with the same mean and SCV}). The second relates to an attempt to substantially reduce the computation times needed, by ignoring a part of the cost function (the so-called {\it value-to-go}).}
  
  \subsection{Non-phase-type service times}
  \tred{
  Concretely, we systematically study the inaccuracy that may be introduced due to replacing the service-time distributions by their phase-type counterparts.  One would like to know how well our method works if the actual service times are not of the phase type, but rather e.g.\ lognormal (motivated by various empirical studies in the healthcare \tredd{context}, such as \cite{CV}) or Weibull distributed.}
  
  \tred{To be able to quantitatively assess this issue, we need to have access to optimal dynamic schedules for the situation that the clients' service times have lognormal or Weibull distributions. We have done this by discretizing time and expressing all quantities relevant in the dynamic programming formulation of the optimal dynamic schedule in terms of convolutions. This technique has also been applied in related papers (cf. \cite{BEG, DEC, VUY, ZY}). Then existing software for dealing with convolutions is used for the numerical evaluation. We refer to \tredd{Appendix \ref{app:alg}} for a compact description of the underlying algorithm. 
  }  
  
  \tred{We performed the following experiment. Define by $\gamma_{\rm XY}$ the total cost when the dynamic schedule is based on the service times having distribution ${\rm X}$, while in reality, they have distribution ${\rm Y}$; here the indices ${\rm X}$ and ${\rm Y}$ are P (i.e., phase type), L (i.e., lognormal), or W (i.e., Weibull). These costs are determined by simulating the system under the dynamic schedule that was determined under ${\rm X}$ using service times that are sampled under ${\rm Y}.$ Table \ref{T7} presents the output (where we have chosen $n = 15$ and $\omega=0.5$).}

 {\small
  \begin{table}[h]
  \centering
\begin{tabular}{|c|c c |cc | c|}
\hline
 \:\:${\mathbb S}(B)$\:\:&\:\:$\gamma_{\rm PL}$ \:\:&\:\: $\gamma_{\rm LL}$\:\: &\:\: $\gamma_{\rm PW}$\:\: &\:\: $\gamma_{\rm WW}$\:\: &\:\: $\gamma_{\rm PP}$\:\: \\ 
 \hline
0.50& 4.16 & 4.15 & 4.38 & 4.38 & 4.34 \\
0.75& 4.98 & 4.96 & 5.31 & 5.31 & 5.31 \\
1.00& 5.62 & 5.60 & 6.05 & 6.05 & 6.05 \\
1.25& 6.14 & 6.12 & 6.66 & 6.66 & 6.57 \\
1.50& 6.58 & 6.57 & 7.20 & 7.19 & 7.03 \\
 \hline 
\end{tabular}
  \caption{\label{T7}Cost of basing optimal dynamic schedule on a phase-type distribution, while the actual service-time distribution is lognormal or Weibull.}
  \end{table}}

\tred{To assess the accuracy of our phase-type based method, one needs to compare the columns $\gamma_{\rm PL}$ and $\gamma_{\rm LL}$ in the lognormal case, and the columns $\gamma_{\rm PW}$ and $\gamma_{\rm WW}$ in the Weibull case. 
The remarkably strong agreement indicates that the dynamic schedules based on the phase-type approximation perform just marginally worse than the dynamic schedules based on the actual distribution (lognormal or Weibull). The numbers in the table provide a strong justification for replacing the service times with their phase-type counterpart. 
We note that the above findings are in line with what was found in the corresponding experiments for the static case as reported in \cite[Section 2.6.1]{AThesis}.}

\tred{In the last column of Table \ref{T7} we also included the cost in case of phase-type service times, with the dynamic schedule being based on the same phase-type distribution. The numbers reveal that in some cases there are substantial differences between $\gamma_{\rm PP}$ and the values in the other columns. The conclusion from this observation is that, while using the phase-type-based optimal schedule leads to different costs in the lognormal, Weibull, and phase-type cases (compare the columns $\gamma_{\rm PL}$, $\gamma_{\rm PW}$, and $\gamma_{\rm PP}$), in all three cases these numbers are close to the respective optimal costs. }

\tred{
The main conclusion from the above experiment is that virtually no accuracy is gained when the dynamic schedules would have been based on the actual service-time distributions rather than their phase-type counterparts. In addition, as pointed out in \cite{SCH}, the alternative method used to obtain the schedules for lognormal and Weibull service times (i.e, discretization of the service times and numerical evaluation of the convolutions) is markedly slower than our approach. As an indication,  for small instances with $n=4$ the computation time of our approach is about 60\% of the computation time of the alternative approach; for $n=8$ this is about 35\%, for $n=12$ about 17\%, and this percentage  decreases even further for larger $n$. In these experiments, we tuned the underlying numerical algorithms such that both approaches reached the same level of precision. It is in addition noted that, when increasing this level of precision, the relative advantage (in terms of computation times) of our method even increased. }

\subsection{Ignoring the value-to-go}\tred{
A naïve, and obviously suboptimal, approach is to schedule the next \tredd{client} without considering future arrivals. More concretely, one could determine the length of the interval until the next arrival, say $x$, based on minimizing $\omega\,f_k(x) + (1-\omega )\,g_k(x)$, but now with $g_k(x)$ defined as the waiting time of the next client; we thus ignore the so-called {\it value-to-go}. The attractive aspect of this approach is the low computational effort needed. Observe that, as the effect of future clients is not incorporated, the dynamic schedule depends on the current number of clients $k$ and the elapsed service time $u$ (of the client in service, that is) only; the total number of clients $n$ does not play a role.}

\tred{We proceed by providing further details on this approach, in which the value-to-go is ignored.
Suppose after a client's arrival there are $k$ clients in the queue, of which there was one client already in service (her elapsed service time being $u$). Let $S\equiv S(k,u)$ denote the {\it total} (remaining) service time of these clients, with $f_S(\cdot)$ its density and $F_S(\cdot)$ its cumulative distribution function; these can be evaluated using the techniques of Section \ref{phase}, as we will explain below. Then the mean idle time until the next client's arrival is
\[{\mathbb I}(x) := {\mathbb E}(x-S)^+ = \int_0^x (x-y)\,f_S(y)\,{\rm d}y,\]
and likewise, mean waiting time of the next client
\[{\mathbb W}(x):= {\mathbb E}(S-x)^+ = \int_x^\infty (y-x)\,f_S(y)\,{\rm d}y.\]
Now performing the optimization, in which we recognize the classical {\it news vendor problem} \cite{STEV},
\[\min_{x\geqslant 0} \omega \,{\mathbb I}(x) + (1-\omega) \,{\mathbb W}(x),\]
we obtain (by differentiating with respect to $x$, and equating the resulting expression to $0$) that the optimizing $x$ solves
\[\omega \,{\mathbb P}(S<x) = (1-\omega) \,{\mathbb P}(S\geqslant x),\]
or, equivalently,
\[x^\star= F_S^{-1}(1-\omega);\]
in words, this means that the optimal interarrival time is the $(1-\omega)$-quantile of $S$. 
This is a highly natural formula, with an appealing interpretation: it quantifies how much the interarrival time should go up when decreasing $\omega$ (evidently, the more we care about idle times, the smaller the interarrival time should be). We thus conclude that for $\omega=\frac{1}{2}$ the optimal interarrival time $x^\star$ equals the {\it median} of $S$; as an aside we mention that, interestingly, it is also argued in \cite{KKM} that if we had wanted to minimize the sum of the {\it second moments} of ${\mathbb I}(x)$ and ${\mathbb W}(x)$ (again with $\omega=\frac{1}{2}$), then $x^\star$ equals the {\it mean} ${\mathbb E}S.$}

\tred{In \cite{KKM} this idea of ignoring the value-to-go has been analyzed in detail for the {\it static} setting; there it is called {\it sequential optimization} (where minimizing the {\it total} cost is called {\it simultaneous optimization}). As pointed out in \cite{KKM}, the sequential problem optimizes an intrinsically different objective function than its simultaneous counterpart. In the static context of \cite{KKM} it is for instance observed that the interarrival times do {\it not} have the dome shape; this makes sense, as the drop at the right-part of the dome shape is essentially due to the fact that there are relatively few remaining \tredd{clients} to be scheduled, i.e., information one does not have in the sequential setup. }

\tred{
The above considerations also mean that in our dynamic context, it is not clear how to rigorously justify the use of a sequential approach, other than its low computational burden. Indeed, as soon as one is able to determine the distribution function of $S\equiv S(k,u)$, for any $k$ and $u$, one can apply this rule, without the need to solve a complex dynamic programming problem. Importantly, with the formulas we have derived, it is straightforward to determine the distribution of $S$. Concretely, it can be checked that $F_{S(k,u)}(x)$ is the $P_{k,u}^\downarrow(x)$ computed in Sections \ref{WED} (for the weighted Erlang case) and \ref{HED} (for the hyperexponential case).}

\tred{
We conclude this subsection by reporting on a series of numerical experiments, quantifying the increase in the total cost when using the sequential approach described above (instead of the optimal dynamic schedule). We performed numerical computations for various values of $n$, ${\mathbb S}(B)$, and $\omega$. The first conclusion is that for small values of $\omega$, the value-to-go can safely be ignored; this makes sense as for these $\omega$ the schedule will be such that the number of waiting clients remains low. In the second place, the effect of ignoring the value-to-go becomes more pronounced with increasing ${\mathbb S}(B)$. And, finally, there is little effect of the number of clients $n$ on the (percentage-wise) increase of the cost. As an illustration, Table \ref{tab:gamma} provides the cost $\gamma^-$ that ignores the value-to-go, the cost $\gamma$ of our approach, and the ratio $\gamma/\gamma^-$, for $n=10$ and $\omega=0.9$.}

 {\small
  \begin{table}[h]
  \centering
\begin{tabular}{|c|c c |c |}
\hline
 \:\:${\mathbb S}(B)$\:\:&\:\:$\gamma^-$ \:\:&\:\: $\gamma$\:\: &\:\: $\gamma/\gamma^-$\:\:  \\ 
 \hline
0.25& 0.88 & 0.87 & 0.99 \\  
0.50& 1.22 & 1.19 & 0.98 \\
0.75& 1.50 & 1.44 & 0.96 \\
1.00& 1.86 & 1.60 & 0.86 \\
1.25& 1.99 & 1.68 & 0.84 \\
1.50& 2.10 & 1.76 & 0.84 \\
1.75& 2.19 & 1.82 & 0.83 \\
2.00& 2.32 & 1.88 & 0.81 \\
 \hline 
\end{tabular}
  \caption{\label{tab:gamma}Cost of optimal dynamic schedule and the variant that ignores the value-to-go, as well as their ratio, for $n=10$ and $\omega=0.9$.}
    \end{table}
    }

\section{Concluding remarks}
In this paper, we have studied dynamic appointment scheduling, focusing on a setting in which at every arrival the arrival epoch of the next client is determined. The fact that the service times may have distributions with an SCV substantially different from 1 has been dealt with by working with phase-type distributions. Therefore, our approach extends beyond the conventional framework in which service times would be exponentially distributed. It led us to resolve the specific complication that the elapsed service time of the client in service has to be accounted for. We rely on dynamic programming to produce the optimal dynamic schedule, where the state information consists of the number of clients in the system and the elapsed service time of the client in service (if any), besides the number of clients still to be scheduled. We have developed an applet that in real time provides the optimal dynamic schedule. Experiments have been performed that provide insight into the achievable gain as a function of the model parameters. 

\tred{The model presented in this paper can be seen as a base model that can be enhanced by adding realistic features. In the first place, borrowing ideas developed in \cite{TEO, KMMB}, no-shows and additional non-anticipated arrivals can be included in our framework. As explained in detail in \cite{KMMB}, this can be done in a natural way by adjusting the parameters of the service times. Overtime can be dealt with similarly, as pointed out in e.g.\ \cite[Section 6.4.3]{AThesis}. Also, it is worthwhile to explore including non-punctuality, potentially building on ideas developed in \cite{KKM, ZY}. }

\textcolor{black}{It is important to notice that in some settings (such as the surgery example and the repairman example mentioned in the introduction), there could be two `stages of waiting’: (1) a `standby mode’ in which the customer has to be at the service location or close to their home, and (2) a true `waiting mode’ in which service is imminent and, for example, anesthesia and drug treatment starts or the customer has to be in the home expecting the repairperson to turn up. The costs of being in the waiting mode for a long time are likely to be much higher than those in the standby mode. While our current model does not explicitly incorporate `standby time' explicitly, in a more advanced model it could be included.}

In our framework, at any arrival epoch, we schedule the arrival epoch of the next client. In the introduction, we provided a set of application areas in which such a mechanism is natural. An evident alternative is to {\it periodically} reschedule, for instance, every hour. It is noted, though, that the latter approach has the intrinsic drawback that at some rescheduling moments potentially {\it no} client is present. Given that it typically takes the client some time to arrive at the service facility, this is an undesirable feature (that the approach that we presented in this paper does not have).

It is anticipated that the framework presented in this paper can be extended in various ways. 
For instance, 
a challenging variant concerns the situation in which at every arrival not the next client, but the one after that is scheduled.
One could in addition think of a setup in which the full remaining schedule is periodically recalculated, for instance, every hour; at every adaptation, the clients are sent an update of their scheduled arrival epochs. In another version, potentially useful for a delivery service, at any adaptation also the order of the clients can be optimized.


\bibliographystyle{plain}

{\small \begin{thebibliography}{111}

\bibitem{AS}{\sc S. Asmussen} (2003). {\it Applied Probability and Queues}, 2nd edition. Springer, New York. 

\bibitem{ANO}
{\sc S. Asmussen, O. Nerman, M. Olsson} (1996). Fitting phase-type distributions via the EM algorithm. {\it Scandinavian Journal of Statistics}, {Vol. 23}, pp. 419-441.

\bibitem{BAI}
{\sc N. Bailey} (1952). A study of queues and appointment systems in hospital out-patient departments, with special reference to waiting-times. {\it Journal of the Royal Statistical Society, Series B (Methodological)}, Vol. 14, pp. 185-199.

\bibitem{BEG}
{\sc M. Begen, M. Queyranne} (2011).
Appointment Scheduling with Discrete Random Durations.
{\it Mathematics of Operational Research,} {Vol. 36}, pp. 240-257.

\bibitem{BELL}
{\sc R. Bellman} (1957).
{\it Dynamic Programming.}
Princeton University Press, New York.

\bibitem{BER}
{\sc B. Berg, B. Denton, S. Erdogan, T. Rohleder, T. Huschka}  (2014). Optimal booking and scheduling in outpatient procedure centers. {\it Computers \& Operations Research}, Vol.  50, pp. 24-37.

\bibitem{CAY}
{\sc T. Cayirli, E. Veral} (2003). Outpatient scheduling in health care: a review of literature. {\it Production and Operations Management}, Vol. 12, pp.\ 519–549.

\bibitem{CV}
{\sc T. Cayirli, E. Veral, H. Rosen} (2006). Designing appointment scheduling systems for
ambulatory care services. {\it  Health Care Management Science}, 
Vol. 9, pp.\ 47-58.

\bibitem{CR}
{\sc R. Chen, L. Robinson}  (2014). Sequencing and scheduling appointments with potential call‐in patients. {\it Production and Operations Management}, {Vol. 23}, pp.\ 1522-1538.

\bibitem{DAV}
{\sc P. Davis} (1972). Gamma Function and Related Functions. In: {\it Handbook of Mathematical Functions with Formulas, Graphs and Mathematical Tables}, Dover Publications, New York, pp. 253-266.

\bibitem{DEC}
{\sc M. Deceuninck, D. Fiems, S. de Vuyst} (2018).
Outpatient scheduling with unpunctual patients and no-shows.
{\it European Journal of Operational Research,} {Vol. 265}, pp. 195-207.

\bibitem{ED}
{\sc S. Erdo\u{g}an, B. Denton} (2013). Dynamic appointment scheduling of a stochastic server with uncertain demand. {\it INFORMS Journal on Computing,} {Vol. 25}, pp.\ 116-132.

\bibitem{GIL}
{\sc J. Gilbertson} (2016).
{\it Queues with a Dynamic Schedule}. Master thesis, University of Melbourne. 

\bibitem{GUP}
{\sc D. Gupta, B. Denton} (2008). Appointment scheduling in health care: Challenges and opportunities. {\it IIE Transactions}, Vol. 40, pp. 800-819.

\bibitem{HG}
{\sc S. Hahn-Goldberg} (2014). {\it Dynamic Optimization Addressing Chemotherapy Outpatient Scheduling}. PhD Thesis, University of Toronto. \url{https://tspace.library.utoronto.ca/bitstream/1807/68108/1/Hahn-Goldberg_Shoshana_201406_PhD_thesis.pdf}

\bibitem{HM}
{\sc R. Hassin, S. Mendel} (2008). Scheduling arrivals to queues: A single-server model with no-shows. {\it Management Science}, Vol.  54, pp. 565-572.

\bibitem{KK}
{\sc G. Kaandorp, G. Koole} (2007). Optimal outpatient appointment scheduling. {\it Health Care Management Science}, Vol. 10, pp. 217-229.

\bibitem{KEMP} {\sc M. de Kemp, M. Mandjes, N. Olver} (2022).
Performance of the smallest-variance-first rule in appointment sequencing. To appear in: {\it Operations Research}. \url{https://pubsonline.informs.org/doi/abs/10.1287/opre.2020.2025}

\bibitem{KKM}
{\sc B. Kemper, C. Klaassen, M. Mandjes} (2014).
Optimized appointment scheduling.
{\it 
European Journal of Operational Research}, Vol.\ {239}, pp.\  243-255.

\bibitem{KONG}
{\sc Q. Kong, C. Lee, C. Teo, Z. Zheng} (2016). Appointment sequencing: Why the smallest-variance-first rule may
not be optimal. {\it European Journal of Operational Research},  Vol. 255, pp. 809-821.

\bibitem{TEO}
{\sc Q. Kong, S. Li, N. Liu, C. Teo, Z. Yan} (2020).
Appointment scheduling under time-dependent patient no-show behavior.
{\it Management Science,} {Vol. 66}, pp. 3480-3500.

\bibitem{AThesis}
{\sc A. Kuiper} (2016).
{\it Optimal Appointment Scheduling in Healthcare.} PhD thesis, University of Amsterdam.  \url{https://pure.uva.nl/ws/files/2776103/174963_AlexKuiper_Thesis_complete.pdf}

\bibitem{KUIP}
{\sc A. Kuiper, B. Kemper, M. Mandjes} (2015).
A computational approach to optimized appointment scheduling.
{\it Queueing Systems}, Vol.\ 79, pp.\ 5-36.

\bibitem{KMM}
{\sc A. Kuiper, M. Mandjes, J. de Mast} (2017).
Optimal stationary appointment schedules.
{\it Operations Research Letters}, 
Vol. 45, pp.\ 549-555.

\bibitem{KMMB} {\sc A. Kuiper, M. Mandjes, J. de Mast, R. Brokkelkamp} (2022). A flexible and optimal approach for appointment scheduling in healthcare. To appear in: {\it  Decision Sciences}. \url{https://onlinelibrary.wiley.com/doi/epdf/10.1111/deci.12517}

\bibitem{MAK}
{\sc H. Mak, Y. Rong, J. Zhang} (2015). Appointment scheduling with limited distributional information. {\it Management Science}, Vol.  61, pp. 316-334.

\bibitem{PEGD}
{\sc C. Pegden, M. Rosenshine} (1990). Scheduling arrivals to
queues. {\it Computers \& Operations Research}, Vol.  17, pp. 343-348.

\bibitem{SCH}
{\sc J. Schwarz, G. Selinka, R. Stolletz} (2016).
Performance analysis of time-dependent queueing systems: Survey
and classification
{\it Omega}, Vol.  63, pp. 170-189.


\bibitem{STEV}
{\sc W. Stevenson} (2013). {\it Operations Management}. McGraw-Hill, New York. 


\bibitem{tijms}
\textsc{H. Tijms} (1994).
\textit{Stochastic Models: An Algorithmic Approach}.
Wiley, New York.



\bibitem{TS}
{\sc M.-Y. Tsang, K. Shehadeh} (2023).
Stochastic optimization models for a home service routing and appointment scheduling problem with random travel and service times. {\it  European Journal of Operational Research} {\bf 307}, 48-63.

\bibitem{VUY}
{\sc S. de Vuyst, H. Bruneel, D. Fiems} (2014).
Computationally efficient evaluation of appointment schedules in health care.
{\it European Journal of Operational Research,} {Vol. 237}, pp. 1142-1154.

\bibitem{WF}
{\sc J. Wang,  R.  Fung} (2015). Dynamic appointment scheduling with patient preferences and choices. {\it Industrial Management \& Data Systems}, Vol.\ 115, pp. 700-717.

\bibitem{W1}
{\sc P. Wang} (1993). Static and dynamic scheduling of customer arrivals to a single-server system. {\it Naval Research Logistics}, Vol.  40, pp. 345-360.

\bibitem{W2}
{\sc P. Wang} (1997). Optimally scheduling $N$ customer arrival times for a single-server system. {\it Computers \& Operations Research}, Vol.  24, pp. 703-716.

\bibitem{WLW}
{\sc S. Wang, N. Liu, G. Wan} (2020). Managing appointment-based services in the presence of walk-in customers. {\it Management Science}, Vol. 66, pp. 667-686.

\bibitem{ZY}
{\sc C. Zacharias, T. Yunes} (2020). Multimodularity in the stochastic appointment scheduling problem with discrete arrival epochs. {\it Management Science,} {Vol. 66}, pp.\ 744-763.

\bibitem{ZAN}
{\sc A. Zander, U. Mohring} (2016).
Dynamic appointment scheduling with patient time preferences and different service time lengths. {\it Lecture Notes in Management Science}, Vol. 8, pp. 72-77.  

\bibitem{ZWW}
{\sc Y. Zhan, Z. Wang, G. Wan} (2021).
Home service routing and appointment scheduling with stochastic service times.
{\it European Journal of Operational Research,} {Vol. 288}, pp. 98-110.

\end{thebibliography}}

\appendix

\section{\tredd{Several proofs}}\label{app:proofs}

\subsubsection*{Proof of Lemma \ref{lem:sum_to_1}}
Instead of the derived recursion for $\varphi_{k\ell}(s)$, the proof makes use of the definition of the generalized Erlang distribution, which immediately yields
\[
c_{k\ell j} = \prod_{\substack{i = k\\ i\neq j}}^{k+\ell}\left(\frac{\mu_i}{\mu_i - \mu_j}\right) \mu_j.
\]
Given our distinct $\mu_k,\dots,\mu_{k+\ell}$ and any function $f$, define the polynomial with coefficients
\[
P(x) := \sum_{j=k}^{k+\ell}L_j(x)f(\mu_j),\quad L_j(x) := \prod_{\substack{i=k\\ i\neq j}}^{k+\ell}\frac{x - \mu_{i}}{\mu_{j} - \mu_{i}}.
\]
It is directly verified that $P(\mu_i) = f(\mu_i)$ for all $i=k,\dots,k+\ell$. If $f$ is a polynomial of degree at most $\ell$ (notation $f\in \mathcal{P}_{\ell}$), then the polynomial is \textit{exact}, i.e., $f = p$. As both $P,f \in \mathcal{P}_{\ell}$, we have $P - f\in \mathcal{P}_{\ell}$. By the fundamental theorem of algebra, it now follows that $P - f = 0$, since this polynomial is of degree at most $\ell$ but has at least $\ell+1$ distinct roots, namely $\mu_k,\dots,\mu_{k+\ell}$. Now take $f(x) = 1$ for all $x$. As $f\in \mathcal{P}_{\ell}$,
\[
1 = P(0) = \sum_{j=k}^{k+\ell}L_j(0)f(0) = \sum_{j=k}^{k+\ell}\prod_{\substack{i=k\\ i\neq j}}^{k+\ell}\frac{\mu_i}{\mu_i - \mu_j},
\]
as desired.

\subsubsection*{Proof of Lemma \ref{lem:3}}
We give a formal proof of the claim, remarking that various other approaches can be followed. By conditioning on the value of ${\rm E}(k,\mu)$, it is directly verified that
\begin{align*}
\psi_{vt}[k,\ell]
&= \int_{t-v}^{t}\mathbb{P}({\rm E}(\ell-k,\mu) > t-x)\,\mathbb{P}({\rm E}(k,\mu) \in {\rm d}x)
= \sum_{i=0}^{\ell-k-1}e^{-\mu t}\frac{\mu^{i+k}}{i!(k-1)!}\int_{t-v}^{t}x^{k-1}(t-x)^i \mathrm{d}x\\
&= \sum_{j=k}^{\ell-1}e^{-\mu t}\frac{\mu^{j}}{(j-k)!(k-1)!}\int_{t-v}^{t}x^{k-1}(t-x)^{j-k} \mathrm{d}x\\
&=\sum_{j=k}^{\ell-1}\mathbb{P}({\rm Pois}(\mu t) = j)\,\frac{j!}{(j-k)!(k-1)!}\int_{1-v/t}^{1}u^{k-1}(1-u)^{j-k} \mathrm{d}u,
\end{align*}
where we used the substitution $j:=i+k$ in the third equality, and the substitution  $u := x/t$
in the last equality. Then recall the property  that
\[
\int_{x}^{1}u^{k-1}(1-u)^{i-1} \mathrm{d}u
= \text{B}(k,i) - \text{B}_{x}(k,i),
\]
where, for $x\in[0,1]$ and $k,i\in{\mathbb N}$,
\[
\text{B}_{x}(k,i) := \int_{0}^{x}t^{k-1}(1 - t)^{i-1}\mathrm{d}t\quad \text{and}
\quad
\text{B}(k,i) := \text{B}_{1}(k,i) = \frac{(k-1)!(i-1)!}{(k+i-1)!},\quad a,b\in \mathbb{N},
\]
are the incomplete beta function and the beta function, respectively. It is a well-known property  \cite{DAV} of these functions that
\[
\text{B}_{x}(k,i) = \text{B}(k,i)I_{x}(k,i) = \text{B}(k,i)\left(1 - I_{1-x}(i,k)\right),
\]
where in particular it holds that, for any $x\in[0,1]$ and $k,i\in\mathbb{N}$,
$
I_x(k,i) = \mathbb{P}(\text{Bin}(k+i-1,x) \geqslant i).
$
Combining the above findings, after some standard algebra, we obtain the stated.
\subsubsection*{Proof of Lemma \ref{L4}}
  By an elementary conditioning argument,
  \begin{align*}\sigma_t[m,k] &= \int_0^t \int_0^s\mu_1e^{-\mu_1 u}\frac{(\mu_1 u)^{m-1}}{(m-1)!}  \left(1-\sum_{\ell=0}^{k-1} e^{-\mu_2(s-u)}\frac{(\mu_2(s-u))^\ell}{\ell!}\right) {\rm d}u\,{\rm d}s\\
  &=
  \int_0^t \int_0^{t-u}\mu_1e^{-\mu_1 u}\frac{(\mu_1 u)^{m-1}}{(m-1)!}  \left(1-\sum_{\ell=0}^{k-1} e^{-\mu_2 s}\frac{(\mu_2 s)^\ell}{\ell!}\right) {\rm d}s\,{\rm d}u\\
  &=  \int_0^t \mu_1e^{-\mu_1 u}\frac{(\mu_1 u)^{m-1}}{(m-1)!}  \left(t-u-\sum_{\ell=0}^{k-1}\int_0^{t-u} e^{-\mu_2 s}\frac{(\mu_2 s)^\ell}{\ell!}{\rm d}s\right) {\rm d}u\\
  &=
t\left(1-\sum_{i=0}^{m-1} e^{-\mu_1 t}\frac{(\mu_1t)^i}{i!}\right)- 
\frac{m}{\mu_1}\left(1-\sum_{i=0}^{m} e^{-\mu_1 t}\frac{(\mu_1t)^i}{i!}\right)\:
-\\ 
&\:\:\:\:\:\:
\sum_{\ell=0}^{k-1} \int_0^t\mu_1e^{-\mu_1 u}\frac{(\mu_1 u)^{m-1}}{(m-1)!}  \int_0^{t-u}
 e^{-\mu_2 s}\frac{(\mu_2 s)^\ell}{\ell!}{\rm d}s\, {\rm d}u 
    ,\end{align*}
using that, for any $\mu>0$, $\ell\in\{0,1,\ldots\}$, and $t\geqslant 0$,
   \begin{align*}\int_0^{t} e^{-\mu s}\frac{(\mu s)^\ell}{\ell!}{\rm d}s&=\frac{1}{\mu}{\mathbb P}\left(
   {\rm E}(\ell+1,\mu)\leqslant t\right)=\frac{1}{\mu}\left(1-\sum_{i=0}^\ell e^{-\mu t}\frac{(\mu t)^i}{i!}\right).
   \end{align*} 
   Using the same argument, it follows that $\sigma_t[m,k]$ equals
     \begin{align*}\lefteqn{
      \left(t-\frac{k}{\mu_2}\right)\left(1-\sum_{i=0}^{m-1} e^{-\mu_1 t}\frac{(\mu_1t)^i}{i!}\right)
      - \frac{m}{\mu_1}\left(1-\sum_{i=0}^{m} e^{-\mu_1 t}\frac{(\mu_1t)^i}{i!}\right)
      \:+}\\
      &\:\:\:\:\:\:\frac{\mu_1}{\mu_2}\sum_{\ell=0}^{k-1}\sum_{i=0}^\ell      
      \int_0^te^{-\mu_1 u}\frac{(\mu_1 u)^{m-1}}{(m-1)!}
      e^{-\mu_2(t-u)}\frac{(\mu_2(t-u))^i}{i!}{\rm d}u\\&= \left(t-\frac{k}{\mu_2}\right)\left(1-\sum_{i=0}^{m-1} e^{-\mu_1 t}\frac{(\mu_1t)^i}{i!}\right)
      - \frac{m}{\mu_1}\left(1-\sum_{i=0}^{m} e^{-\mu_1 t}\frac{(\mu_1t)^i}{i!}\right)+\frac{\mu_1}{\mu_2}\sum_{\ell=0}^{k-1}\sum_{i=0}^\ell      
       \rho_t[m-1,i].     
     \end{align*}
  This proves the claim.

\subsubsection*{Proof of Lemma \ref{L5}}
 Observe that, by distinguishing between $s\in[0,t-v)$ and $s\in[t-v,t]$,
 \begin{align*}
 \chi_{vt,1}[k,\ell] &=  \int_0^\infty{\mathbb P}(t-v -s\leqslant {\rm E}(k,\mu_1)\leqslant t-s,
 {\rm E}(k,\mu_1)+ {\rm E}(1,\mu_1)>t-s) \,{\mathbb P}({\rm E}(\ell,\mu_2)\in {\rm d}s)\\
 &=\int_0^{t-v} e^{-\mu_1(t-s)}\frac{(\mu_1(t-s))^k}{k!}{\mathbb P}\left({\rm Bin}\left(k, \frac{v}{t-s} \right)>0\right)
\,{\mathbb P}({\rm E}(\ell,\mu_2)\in {\rm d}s)\:+\\
 &\:\:\:\:\:\:\:\:\int_{t-v}^t  e^{-\mu_1(t-s)}\frac{(\mu_1(t-s))^k}{k!} \,{\mathbb P}({\rm E}(\ell,\mu_2)\in {\rm d}s),
 \end{align*}
 where the first term in the last expression is due to Lemma \ref{lem:3}. 
 Inserting 
\[{\mathbb P}\left({\rm Bin}\left(k, \frac{v}{t-s} \right)>0\right) = 1 - \left(\frac{t-s-v}{t-s}\right)^k,\] 
we find that $ \chi_{vt,1}[k,\ell]$ equals 
\begin{align*}
\mu_2\int_{0}^t  e^{-\mu_1(t-s)}\frac{(\mu_1(t-s))^k}{k!} &\,e^{-\mu_2 s}\frac{(\mu_2 s)^{\ell-1}}{(\ell-1)!} {\rm d}s\\
-&\mu_2\int_0^{t-v} e^{-\mu_1(t-s)}\frac{(\mu_1(t-s))^k}{k!} \left(\frac{t-s-v}{t-s}\right)^k \,e^{-\mu_2 s}\frac{(\mu_2 s)^{\ell-1}}{(\ell-1)!} {\rm d}s.
\end{align*}
The first term we recognize as $\mu_2\,\rho_{t}[k,\ell-1].$ The second term can be rewritten as
\[e^{-\mu_1 v}\mu_2\int_{0}^{t-v} e^{-\mu_1(t-s-v)}\frac{(\mu_1(t-s-v))^k}{k!} \,e^{-\mu_2 s}\frac{(\mu_2 s)^{\ell-1}}{(\ell-1)!} {\rm d}s,\]
which equals $e^{-\mu_1 v}\mu_2 \,\rho_{t-v}[k,\ell-1]$.

\section{Numerical Results for the Heterogeneous Exponential Case}
\label{app:hetexp} 

\tredd{
In this appendix, we present some numerical results for the dynamic appointment scheduling model with heterogeneous exponentially distributed service times, as discussed in Section~\ref{hetexp}. 
With heterogeneous service times, there is an additional, interesting challenge to find the optimal scheduling order for the clients. 
As backed by the results in \cite{KEMP}, the so-called {\it smallest-variance-first} rule (ordering the clients such that their mean service times, and hence also the corresponding variances, are increasing) is a logical candidate.
}

\begin{example}\label{E3}{\em 
As in Example \ref{E1}, we wish to quantify the gain achieved by scheduling dynamically rather than in advance. As before, we let $K_{\text{pre}}(n,\omega)$ and $K_{\text{dyn}}(n,\omega)$ be the value of the objective function of the precalculated and  the dynamic schedule, respectively. Again $r(n,\omega)$ denotes the ratio of $K_{\text{dyn}}(n,\omega)$ and $K_{\text{pre}}(n,\omega)$. First of all, we take $n$ equally spaced parameters in the interval $[0.5,1.5]$, i.e., for $i=1,\dots,n$,
\[
\mu_i = 0.5 + \frac{i - 1}{n - 1}.
\]
The results are shown in Table~\ref{het_adap_vs_presc}. One of the conclusions from this table is that the cost grows superlinearly in the number of clients $n$. As in the homogeneous case, the gain is more pronounced when $\omega$ and/or $n$ are large.

{\small
\begin{table}
\centering
\begin{tabular}{|c|c|ccccccccc|}
\hline
$n$   & $\omega$        & 0.1  & 0.2  & 0.3  & 0.4  & 0.5  & 0.6  & 0.7  & 0.8  & 0.9 \\
\hline
5     & $K_{\text{dyn}}(n,\omega)$ & 1.23 & 1.79 & 2.10 & 2.24 & 2.24 & 2.12 & 1.87 & 1.47 & 0.90 \\
      & $K_{\text{pre}}(n,\omega)$ & 1.32 & 2.01 & 2.43 & 2.65 & 2.70 & 2.58 & 2.28 & 1.79 & 1.07 \\
      & $r(n,\omega)$             & 0.93 & 0.89 & 0.87 & 0.85 & 0.83 & 0.82 & 0.82 & 0.82 & 0.85 \\
\hline
10    & $K_{\text{dyn}}(n,\omega)$ & 2.52 & 3.68 & 4.33 & 4.63 & 4.65 & 4.42 & 3.94 & 3.16 & 2.00 \\
      & $K_{\text{pre}}(n,\omega)$ & 2.71 & 4.16 & 5.13 & 5.73 & 6.00 & 5.94 & 5.51 & 4.60 & 3.01 \\
      & $r(n,\omega)$             & 0.93 & 0.88 & 0.84 & 0.81 & 0.78 & 0.74 & 0.71 & 0.69 & 0.67 \\
\hline
15    & $K_{\text{dyn}}(n,\omega)$ & 3.83 & 5.58 & 6.56 & 7.02 & 7.06 & 6.72 & 5.99 & 4.83 & 3.08 \\
      & $K_{\text{pre}}(n,\omega)$ & 4.09 & 6.29 & 7.77 & 8.74 & 9.23 & 9.25 & 8.73 & 7.50 & 5.15 \\
      & $r(n,\omega)$             & 0.94 & 0.89 & 0.84 & 0.80 & 0.77 & 0.73 & 0.69 & 0.64 & 0.60 \\
\hline
\end{tabular}
\caption{Cost of dynamic and precalculated schedule for Example~\ref{E3}, $\mu\in [0.5,1.5]$.
\label{het_adap_vs_presc}}
\end{table}}

\noindent  We proceed by studying the effect of the spread of the parameters. To this end, we take $n=10$ and equally spaced parameters in an interval of length $\Delta$ that is centered around 1. More formally, we choose, for $i=1,\dots,n=10$ and $ \Delta \in (0,2)$,
\[
\mu_i = 1 - \frac{\Delta}{2} + \frac{i-1}{n-1}\Delta.
\]
The obtained results are given in Table \ref{het_delta}. For small values of $\Delta$, the parameters are relatively homogeneous, so that the gains resemble those achieved in the homogeneous case. The main conclusion is that the higher the value of $\Delta$, the more spread in the parameters, the higher the gain.  

{\small
\begin{table}
\centering
\begin{tabular}{|c|c|ccccccccc|}
\hline
$\Delta$   & $\omega$        & 0.1  & 0.2  & 0.3  & 0.4  & 0.5  & 0.6  & 0.7  & 0.8  & 0.9 \\
\hline
0.25  & $K_{\text{dyn}}(\omega)$ & 2.18 & 3.16 & 3.70 & 3.95 & 3.95 & 3.74 & 3.31 & 2.63 & 1.65 \\
      & $K_{\text{pre}}(\omega)$ & 2.30 & 3.49 & 4.25 & 4.69 & 4.87 & 4.77 & 4.37 & 3.61 & 2.33 \\
      & $r(\omega)$             & 0.95 & 0.91 & 0.87 & 0.84 & 0.81 & 0.78 & 0.76 & 0.73 & 0.71 \\
\hline
0.5   & $K_{\text{dyn}}(\omega)$ & 2.25 & 3.27 & 3.83 & 4.09 & 4.10 & 3.88 & 3.44 & 2.75 & 1.73 \\
      & $K_{\text{pre}}(\omega)$ & 2.39 & 3.63 & 4.44 & 4.92 & 5.12 & 5.03 & 4.62 & 3.83 & 2.48 \\
      & $r(\omega)$             & 0.94 & 0.90 & 0.86 & 0.83 & 0.80 & 0.77 & 0.74 & 0.72 & 0.70 \\
\hline
0.75  & $K_{\text{dyn}}(\omega)$ & 2.36 & 3.43 & 4.03 & 4.31 & 4.32 & 4.10 & 3.64 & 2.92 & 1.84 \\
      & $K_{\text{pre}}(\omega)$ & 2.52 & 3.84 & 4.71 & 5.24 & 5.47 & 5.40 & 4.98 & 4.14 & 2.70 \\
      & $r(\omega)$             & 0.94 & 0.89 & 0.86 & 0.82 & 0.79 & 0.76 & 0.73 & 0.70 & 0.68 \\
\hline
1     & $K_{\text{dyn}}(\omega)$ & 2.52 & 3.68 & 4.33 & 4.63 & 4.65 & 4.42 & 3.94 & 3.16 & 2.00 \\
      & $K_{\text{pre}}(\omega)$ & 2.71 & 4.16 & 5.13 & 5.73 & 6.00 & 5.94 & 5.51 & 4.60 & 3.01 \\
      & $r(\omega)$             & 0.93 & 0.88 & 0.84 & 0.81 & 0.78 & 0.74 & 0.71 & 0.69 & 0.67 \\
\hline
1.25  & $K_{\text{dyn}}(\omega)$ & 2.79 & 4.07 & 4.79 & 5.14 & 5.18 & 4.93 & 4.40 & 3.55 & 2.26 \\
      & $K_{\text{pre}}(\omega)$ & 3.02 & 4.68 & 5.79 & 6.51 & 6.85 & 6.82 & 6.35 & 5.33 & 3.51 \\
      & $r(\omega)$             & 0.92 & 0.87 & 0.83 & 0.79 & 0.76 & 0.72 & 0.69 & 0.67 & 0.64 \\
\hline
1.5   & $K_{\text{dyn}}(\omega)$ & 3.26 & 4.78 & 5.65 & 6.07 & 6.13 & 5.86 & 5.25 & 4.26 & 2.73 \\
      & $K_{\text{pre}}(\omega)$ & 3.61 & 5.65 & 7.06 & 7.99 & 8.46 & 8.47 & 7.93 & 6.70 & 4.43 \\
      & $r(\omega)$             & 0.90 & 0.85 & 0.80 & 0.76 & 0.72 & 0.69 & 0.66 & 0.64 & 0.62 \\
\hline
\end{tabular}
\caption{Cost of dynamic and precalculated schedule for Example~\ref{E3}, $\mu \in [1-\frac{\Delta}{2},1+\frac{\Delta}{2}]$.
\label{het_delta}}
\end{table}}

As announced before, we can now also analyze the effect of the order in which the clients are scheduled. We do so for equally spaced $\mu_i$ in the interval $[0.5,1.5]$ and $n=10$. We first compute the cost when the $\mu_i$ are increasing, corresponding to the smallest-variance-first strategy studied in \cite{KEMP}. Secondly, we pick a randomly selected permutation of the $\mu_i$. Finally, we take decreasing $\mu_i$. The results are given in Table \ref{het_smallestvar}. The main conclusion is that the numerics align with findings of \cite{KEMP}, in the sense that the cost of serving clients with decreasing $\mu_i$ is substantially lower than for increasing $\mu_i$ or randomly ordered $\mu_i.$

{\small \begin{table}
\centering
\begin{tabular}{|c|c|ccccccccc|}
\hline
$\mu$ & $\omega$        & 0.1  & 0.2  & 0.3  & 0.4  & 0.5  & 0.6  & 0.7  & 0.8  & 0.9 \\
\hline
increasing & $K_{\text{dyn}}(\omega)$ & 2.52 & 3.68 & 4.33 & 4.63 & 4.65 & 4.42 & 3.94 & 3.16 & 2.00 \\
           & $K_{\text{pre}}(\omega)$ & 2.71 & 4.16 & 5.13 & 5.73 & 6.00 & 5.94 & 5.51 & 4.60 & 3.01 \\
           & $r(\omega)$             & 0.93 & 0.88 & 0.84 & 0.81 & 0.78 & 0.74 & 0.71 & 0.69 & 0.67 \\
\hline
random     & $K_{\text{dyn}}(\omega)$ & 2.30 & 3.35 & 3.94 & 4.21 & 4.23 & 4.02 & 3.57 & 2.86 & 1.80 \\
           & $K_{\text{pre}}(\omega)$ & 2.49 & 3.82 & 4.68 & 5.20 & 5.42 & 5.33 & 4.91 & 4.06 & 2.62 \\
           & $r(\omega)$             & 0.92 & 0.88 & 0.84 & 0.81 & 0.78 & 0.75 & 0.73 & 0.70 & 0.69 \\
\hline
decreasing & $K_{\text{dyn}}(\omega)$ & 2.18 & 3.15 & 3.67 & 3.38 & 3.86 & 3.62 & 3.17 & 2.49 & 1.53 \\
           & $K_{\text{pre}}(\omega)$ & 2.27 & 3.38 & 4.06 & 4.42 & 4.51 & 4.35 & 3.92 & 3.17 & 1.99 \\
           & $r(\omega)$             & 0.96 & 0.93 & 0.90 & 0.88 & 0.85 & 0.83 & 0.81 & 0.79 & 0.77 \\
\hline
\end{tabular}
\caption{Cost of dynamic and precalculated schedule  for Example~\ref{E3}, $\mu\in [0.5,1.5]$.}
\label{het_smallestvar}
\end{table}}

\tredd{
It is interesting to conduct a similar experiment as Experiment 1 in Section~\ref{subsec:experiments}, where we have studied the impact of  randomized parameter values on (the quality of) the schedules. This time, we fix the number of clients at $n=10$ and sample the $\mu_i$ values uniformly at random from the interval $[0.5, 1.5]$. 
We compute 95\% confidence intervals (CIs) and prediction intervals~(PIs) for the difference of the costs of precalculated and dynamic scheduling, for $\omega\in\{0.2, 0.5, 0.8\}$, based on 100,000 independent simulation runs.
The results, shown in Table~\ref{tbl:randomparametersheterogeneous}, are comparable to those in Experiment 1, confirming once again that dynamic scheduling strongly decreases the costs of a schedule. The absolute decrease remains approximately constant for $\omega>0.5$, something that can also be observed in Tables~\ref{het_delta} and \ref{het_smallestvar}. 
}
$\hfill\Diamond$

\begin{table}[ht]
\centering
\begin{tabular}{|c|ccc|cc|}
\hline
 $\omega$  &  $K_{\text{dyn}}(\omega)$ & $K_{\text{pre}}(\omega)$ &
 $K_{\text{pre}}(\omega) - K_{\text{dyn}}(\omega)$ &
 95\% CI & 95\% PI\\
\hline      
0.2   & 3.422 & 3.811 & 0.389 & [0.388, 0.391] & [0.268, 0.506] \\
0.5   & 4.279 & 5.273 & 0.994 & [0.993, 0.995] & [0.721, 1.302] \\ 
0.8   & 2.849 & 3.855 & 1.006 & [1.005, 1.007] & [0.733, 1.334]\\
\hline
\end{tabular}
\caption{Cost of dynamic and precalculated schedule for Example~\ref{E3}, with $n=10$ and random $\mu_i$.  95\% confidence intervals and prediction intervals are computed for the difference in costs.}
\label{tbl:randomparametersheterogeneous}
\end{table}

}\end{example}

\section{General Case: Mean Idle Time}\label{app:mean_idle_time}

\tredd{This appendix contains the computation of the mean idle time $\bar{g}^{\circ}_{k,u}(t)$ in the general case, i.e., assuming the service times stem from the weighted Erlang and hyperexponential distribution described in Section \ref{sec:phase_type_fit}. For completeness, we formulate  the  dynamic programming algorithm that uses  the $\bar{g}^{\circ}_{k,u}(t)$ in Theorem~\ref{thm:thm5}.}

 \begin{theorem}\label{thm:thm5} We can determine the $C_i(k,u)$ recursively: for $i = 1,\ldots,n-1$, $k = 1,\ldots,i$, and $u\geqslant 0$,
 \begin{align*}C_i(k,u) 
 = \inf_{t\geqslant 0} \Bigg(\omega \,\bar f^\circ_{k,u}(t) &+ (1-\omega) \,\bar g^\circ_{k,u}(t)
 +\sum_{\ell=2}^{k} \int_{(0,t)} p_{k\ell,uv}(t) \,C_{i+1}(\ell,v) \,{\rm d}v\\
 &+\,P^\downarrow_{k,u}(t)\, C_{i+1}(1,0) + P^\uparrow_{k,u}(t) \,C_{i+1}(k+1, u+t)\Bigg),
 \end{align*}
 whereas, for $k=1,\ldots,n$ and $u\geqslant 0$,
 \[C_n(k,u) = (1-\omega) \,\bar g^\circ_{k,u}(\infty).\]
 \end{theorem}

\tredd{The expressions for $\bar{g}^{\circ}_{k,u}(t)$ are derived below, first for the case of weighted Erlang distributed service times (with SCV less than 1) and then for the case of hyperexponential distributed service times (with SCV greater than 1).}

\subsection{Weighted Erlang distribution}
We start by evaluating\tredd{, for $z = 1,\dots,K+1$,}
\begin{align*}\bar g_{k z}(t)&:= \int_0^t\sum_{\ell=0}^{k-1} (k-\ell-1) {\mathbb P}\big(N_s=k-\ell\,|\, N_{0+}=k, Z_{0+}=z\big) {\rm d}s.
\end{align*}
To this end, we first focus on computing $ {\mathbb P}(N_s=k-\ell\,|\, N_{0+}=k, Z_{0+}=z)$. Again we need to distinguish between $z=1,\ldots, K$ and $z=K+1$. The event under consideration corresponds to $\ell$ clients having left at time $s$, but the $(\ell+1)$-st still being in the system, with the client in service at time $0$ being in phase $z$. 

$\circ$~If $z=1,\ldots,K$, then, using the same argument as used when computing $\bar f_{kz}(t)$, for $\ell=1,\ldots,k-1$,
\begin{align}\nonumber {\mathbb P}(N_s=k-\ell\,|\, N_{0+}=k, &\,Z_{0+}=z) = 
\sum_{m=0}^\ell \mathbb{P}({\rm Bin}(\ell,1-p) = m)
{\mathbb P}\big({\rm E}(\ell K- z +1 + m,\mu) \leqslant s\big)\:-\\
&\hspace{-3mm}\sum_{m=0}^{\ell+1} \mathbb{P}({\rm Bin}(\ell+1,1-p) = m) 
 {\mathbb P}\big({\rm E}((\ell+1)K- z +1 + m,\mu) \leqslant s\big).\label{casemain}
\end{align}
For $\ell=0$, we should have that the client in service does not leave before $s$, so that
\begin{align*}{\mathbb P}(N_s=k\,|\, N_{0+}=k, Z_{0+}=z) &= p\, {\mathbb P}\big({\rm E}(K- z +1,\mu) > s\big) +(1-p)\,
 {\mathbb P}\big({\rm E}(K- z +2,\mu) > s\big)\\
&=1- p\, {\mathbb P}\big({\rm E}(K- z +1,\mu) \leqslant s\big) -(1-p)\,
 {\mathbb P}\big({\rm E}(K- z +2,\mu) \leqslant s\big),\end{align*}
 which is in line with \eqref{casemain} for $\ell=0.$ 
It now directly follows that, for $z=1,\ldots,K$,
\begin{align*}
\bar g_{kz}(t) &=  \sum_{\ell=0}^{k-1} (k-\ell-1)\sum_{m=0}^\ell \mathbb{P}({\rm Bin}(\ell,1-p) = m) f_{\ell K- z +1 + m}(t)\:-\\ &\hspace{3mm} \sum_{\ell=0}^{k-1} (k-\ell-1)\sum_{m=0}^{\ell+1} \mathbb{P}({\rm Bin}(\ell+1,1-p) = m)  f_{(\ell+1) K- z +1 + m}(t).
\end{align*}
$\circ$~We proceed by analyzing the case $z=K+1$, along the same lines. As before, we can ignore $\ell=k$, whereas for $\ell=1,\ldots,k-1$,
\begin{align}\nonumber {\mathbb P}(N_s=k-\ell\,|\, &N_{0+}=k,Z_{0+}=z) \\&= 
\sum_{m=0}^\ell \mathbb{P}({\rm Bin}(\ell-1,1-p) = m)
{\mathbb P}\big({\rm E}((\ell-1) K- z +1 + m,\mu) \leqslant s\big)\:-\nonumber\\
&\hspace{5mm}\sum_{m=0}^{\ell+1} \mathbb{P}({\rm Bin}(\ell,1-p) = m) 
 {\mathbb P}\big({\rm E}(\ell K- z +1 + m,\mu) \leqslant s\big).\label{casemain2}\end{align}
The remaining case is $\ell=0$, for which we have
\begin{align*}{\mathbb P}(N_s=k\,|\, N_{0+}=k, Z_{0+}=z) &= {\mathbb P}\big({\rm E}(1,\mu) > s\big) ,\end{align*}
which is {\it not} in line with \eqref{casemain2} for $\ell=0.$ Upon combining the above, we obtain
 \begin{align*}
\bar g_{kz}(t) &=\sum_{\ell=1}^{k-1} (k-\ell-1) \sum_{m=0}^{\ell-1} \mathbb{P}({\rm Bin}(\ell-1,1-p) = m)  f_{(\ell-1) K +1 + m}(t)\:-\\ &\hspace{3mm} \sum_{\ell=1}^{k-1} (k-\ell-1)\sum_{m=0}^{\ell} \mathbb{P}({\rm Bin}(\ell,1-p) = m) f_{\ell K+1 + m}(t)+(k-1)\,\frac{1-e^{-\mu t}}{\mu}.
\end{align*}
In line with the way we translated expressions for $\bar f_{k z}(t)$ into expressions for $\bar f^\circ_{k, u}(t)$, we obtain that
\begin{align*}\bar g^\circ_{k,u}(t)&:= \int_0^t\sum_{\ell=0}^{k-1} (k-\ell-1) {\mathbb P}\big(N_s=k-\ell\,|\, N_{0+}=k, \bar{B}_{0+}=u\big) {\rm d}s= \sum_{z=1}^{K+1} \gamma_z(u) \bar g_{kz}(t).
\end{align*}

\subsection{Hyperexponential distribution}

\tredd{
Again we start by evaluating, for $z = 1,2$,
\begin{align*}\bar g_{k z}(t)&:= \int_0^t\sum_{\ell=0}^{k-1} (k-\ell-1) {\mathbb P}\big(N_s=k-\ell\,|\, N_{0+}=k, Z_{0+}=z\big) {\rm d}s.
\end{align*}
As in the weighted Erlang case, we first compute $ {\mathbb P}(N_s=k-\ell\,|\, N_{0+}=k, Z_{0+}=z)$, corresponding to $\ell$ clients having left at time $s$, but the $(\ell+1)$-st still in the system, given that the service time of the client in service is exponentially distributed with parameter $\mu_z$. For $\ell = 0$ we have
\begin{align*}
{\mathbb P}(N_s=k\,|\, N_{0+}=k, &\,Z_{0+}=z)
= {\mathbb P}\big({\rm E}(1,\mu_z) > s\big),
\end{align*}
whereas for $\ell = 1,\dots,k-1$ we obtain
\begin{align}\nonumber {\mathbb P}(N_s=k-\ell\,|\, N_{0+}=k, &\,Z_{0+}=z) \\ 
&\nonumber\hspace{-22mm}= 
\sum_{m=0}^{\ell-1} \mathbb{P}({\rm Bin}(\ell-1,p) = m)
{\mathbb P}\big({\rm E}(m+{\mathbbm 1}_{\{z=1\}},\mu_1) + {\rm E}(\ell-1-m+{\mathbbm 1}_{\{z=2\}},\mu_2) \leqslant s\big)\:-\\
&\nonumber\hspace{-18mm}\sum_{m=0}^{\ell} \mathbb{P}({\rm Bin}(\ell,p) = m) 
 {\mathbb P}\big({\rm E}(m+{\mathbbm 1}_{\{z=1\}},\mu_1) + {\rm E}(\ell-m+{\mathbbm 1}_{\{z=2\}},\mu_2) \leqslant s\big).
\end{align}
Combining the above yields
\begin{align*}
\bar{g}_{kz}(t) &= 
\sum_{\ell=1}^{k-1} (k-\ell-1)\Bigg(
\sum_{m=0}^{\ell-1} \mathbb{P}({\rm Bin}(\ell-1,p) = m)\sigma_t\big[m+{\mathbbm 1}_{\{z=1\}},\ell-1-m+{\mathbbm 1}_{\{z=2\}}\big] - \\
&\quad \sum_{m=0}^{\ell}
\mathbb{P}({\rm Bin}(\ell,p) = m)\sigma_t\big[m+{\mathbbm 1}_{\{z=1\}},\ell-m+{\mathbbm 1}_{\{z=2\}}\big]
\Bigg)
+ (k-1)\frac{1 - e^{-\mu_z t}}{\mu_z}.
\end{align*}
With the $\bar{g}_{kz}(t)$ readily available, it follows that
\begin{align*}\bar g^\circ_{k,u}(t)&:= \int_0^t\sum_{\ell=0}^{k-1} (k-\ell-1) {\mathbb P}\big(N_s=k-\ell\,|\, N_{0+}=k, \bar{B}_{0+}=u\big) {\rm d}s= \gamma_1(u) \bar g_{k1}(t)
+ \gamma_2(u) \bar g_{k2}(t).
\end{align*}}


\section{Description of convolution-based algorithm}\label{app:alg}

\tred{We here provide a brief account of a method by which one can numerically evaluate a discrete-time approximation of the objective function. In this setup an equidistant  discrete grid is chosen on which the service times are defined; in the description below we assume, without loss of generality,  that the service times attain positive integer values.}

\tred{Let as before $N_t$ denote the number of clients present at time $t\in{\mathbb N}$, and $R_t$ the elapsed service time of the client in service at $t$ (if any). In order to compute the objective function,  
we have to compute probabilities of the type 
\begin{equation}\label{prob0}
\mathfrak{p}[(m,k)\to(n,\ell)]:=
{\mathbb P}(N_t = n , R_t= \ell\,|\,N_0=m, R_0 = k),\end{equation}
and various related quantities (of which some are considerably more involved). 
To compute the probability \eqref{prob0}, we define, for a given $k$, the distribution
\[{\mathbb P}(\xi^\star_k = k') = \frac{{\mathbb P}(\xi= k+k')}{{\mathbb P}(\xi\geqslant k)},\]
where $\xi$ is distributed as a single service time;
the random variable $\xi^\star_k$ represents the residual service time, conditional on the service time being at least~$k$. As is seen easily,  the evaluation of $\mathfrak{p}[(m,k)\to(n,\ell)]$ requires a routine to compute probabilities of the type
\begin{equation}\label{prob}{\mathbb P}\left(\xi^\star_k +\sum_{i=1}^j\xi_i\leqslant t, \xi^\star_k +\sum_{i=1}^{j+1}\xi_i> t \right)
={\mathbb P}\left(\xi^\star_k +\sum_{i=1}^j\xi_i\leqslant t \right)-{\mathbb P}\left( \xi^\star_k +\sum_{i=1}^{j+1}\xi_i\leqslant t \right)
,\end{equation}
where the $\xi_i$ are i.i.d.\ copies of $\xi$, also independent of $\xi_k^\star$.}

\tred{Suppose now that the $\xi_i$ are, say, (the discrete-time version of) lognormal or Weibull, then a complication is that we do not have explicit expressions for the density or the cumulative distribution function of the random variables appearing in the probability \eqref{prob}. This means that we have to resort to numerical techniques to evaluate it, requiring the computation of a $(j+1)$-fold and a $(j+2)$-fold convolution. This is most efficiently done by relying on fast Fourier transform methodology (and that is also how we implemented it).
}

\end{document}